\documentclass[10pt]{amsart}
\usepackage{amsfonts}
\usepackage{graphicx}
\usepackage{amscd}
\usepackage{amsmath,amsfonts,amssymb,amsthm,mathrsfs,xypic}
\xyoption{curve}
\usepackage{fancybox}
\newcommand{\Hom}{\operatorname{Hom}}

\newcommand{\s}{\hfill \blacksquare}
\newtheorem{thm}{Theorem}[section]

\newtheorem{cor}[thm]{Corollary}
\newtheorem{exm}[thm]{Example}

\newtheorem{lem}[thm]{Lemma}
\newtheorem{rem}[thm]{Remark}
\newtheorem{defn}[thm]{Definition}
\newtheorem{prop}[thm]{Proposition}
\newtheorem{exm-prop}[thm]{Example-Proposition}
\begin{document}
\title[Types of  Serre subcategories]
{Types of  Serre subcategories of Grothendieck categories}
\author[Jian FENG, Pu Zhang]{Jian FENG, Pu Zhang$^*$}
\vspace{-10pt}\begin{abstract} \ Every Serre subcategory of an
abelian category is assigned a unique type. The type of a Serre
subcategory of a Grothendieck category is in the list:
$$(0, 0), \ (0, -1), \ (1, -1), \ (0, -2), \ (1, -2), \ (2, -1), \ (+\infty,
-\infty);$$ and for each $(m, -n)$ in this list, there exists a
Serre subcategory such that its type is $(m, -n)$. This uses right
(left) recollements of abelian categories, Tachikawa-Ohtake [TO] on
strongly hereditary torsion pairs, and Geigle-Lenzing [GL] on
localizing subcategories. If all the functors in a recollement of
abelian categories are exact, then the recollement splits. Quite
surprising, any left recollement of a Grothendieck category can be
extended to a recollement; but this is not true for a right
recollement. Thus, a colocalizing subcategory of a Grothendieck
category is localizing; but the converse is not true. All these
results do not hold in triangulated categories. \vskip5pt {\it Key
words and phrases. \ the type of a Serre subcategory, right
recollement, strongly hereditary torsion pair, quotient functor,
localizing subcategory, Grothendieck category}
\end{abstract}
\thanks{$^*$The corresponding author. \ \ fengjian008@sina.com $($J. Feng$)$, \  pzhang@sjtu.edu.cn $($P. Zhang$)$.}
\thanks{Supported by the NSFC 11271251 and 11431010.}
\maketitle
\vspace{-25pt}
\section{\bf Introduction}

Given a Serre subcategory $\mathcal S$ of an abelian category
$\mathcal A$ with inclusion functor $i:\mathcal S \rightarrow
\mathcal A$ and quotient functor $Q:\mathcal A\rightarrow \mathcal
A/\mathcal S$, it is fundamental to know when it is localizing
(resp. colocalizing), i.e., $Q$ has a right (resp. left) adjoint
([S], [G]). By W. Geigle and H. Lenzing [GL], $\mathcal S$ is
localizing if and only if there exists an exact sequence
$0\rightarrow S_1 \rightarrow A \rightarrow C\rightarrow
S_2\rightarrow 0$ with $S_1\in \mathcal S, \ S_2\in \mathcal S,$ and
$C\in \mathcal S^{\perp_{\le 1}};$ and if and only if the
restriction $Q:\mathcal S^{\perp_{\le 1}} \rightarrow \mathcal
A/\mathcal S$ is an equivalence of categories. In this case the
right adjoint of $Q$ is fully faithful. There is a corresponding
work for a thick triangulated subcategory of a triangulated category
(A. Neeman  [N, Chap. 9]).

\vskip5pt

It is then natural to describe Serre subcategories of a fixed
abelian category via the length of two adjoint sequences where $i$
and $Q$ lie. A finite or an infinite  sequence \ $(\cdots, F_{1},
F_0, F_{-1}, \cdots)$ \ of functors between additive categories is
{\it an adjoint sequence}, if each pair $(F_i, F_{i-1})$ is an
adjoint pair. Each functor in an adjoint sequence is additive.

\vskip5pt

Let $\mathcal S$ be a Serre subcategory of an abelian category
$\mathcal A$ with the inclusion functor $i: \mathcal S \rightarrow
\mathcal A$ and the quotient functor $Q: \mathcal A\rightarrow
\mathcal A/ \mathcal S$. The pair $(\mathcal S, i)$ is
of {\it type $(m ,-n)$}, or in short, the Serre subcategory
$\mathcal S$ is of {\it type $(m ,-n)$}, where $m$ and $n$ are in
the set $\Bbb N_0\cup\{+\infty\}$, and $\Bbb N_0$ is the set of
non-negative integers, provided that there exist adjoint sequences
$$(F_m, \cdots, F_{1}, \ F_0 = i, \ F_{-1}, \cdots, F_{-n}) \ \ \ \
\mbox{ and} \ \ \ \ (G_m, \cdots, G_{1}, \ G_0 = Q, \ G_{-1},
\cdots, G_{-n})$$ such that $F_m$ and $G_m$ can not have left
adjoints at the same time, and that $F_{-n}$ and $G_{-n}$ can not
have right adjoints at the same time.

\vskip5pt

We stress that the type of $\mathcal S$ depends on the abelian
category $\mathcal A$ in which $\mathcal S$ is a Serre subcategory.
Since in an adjoint pair one functor uniquely determines the other,
every Serre subcategory is of a unique type $(m, -n)$. We will see a
Serre subcategory $\mathcal S$ of type $(1, -2)$, but with adjoint
sequences \begin{align*}(F_1, \ & \ i, \ \ F_{-1}, \ F_{-2}, \
F_{-3}, F_{-4}, \  F_{-5}), \\ (G_4, \ G_3, \ G_2, \ G_1, \ & Q, \ \
G_{-1}, \ G_{-2}).\end{align*} See Remark \ref{strange}.

\vskip5pt

{\it A right recollement} $(\mathcal B, \ \mathcal A, \ \mathcal C,
\ i_*, \ i^!, \ j^*,\ j_*)$  (see e.g. [P], [K\"o]) of  abelian
categories is a diagram

\begin{center}
\begin{picture}(100,28)
\put(13,10){\makebox(-22,2) {$\mathcal B$}}
\put(10,16){\vector(1,0){30}} \put(40,6){\vector(-1,0){30}}

\put(37,10){\makebox(25,0.8) {$\mathcal A$}}
\put(60,16){\vector(1,0){30}} \put(90,6){\vector(-1,0){30}}
\put(88,10){\makebox(25,0.5){$\mathcal C$}}
\put(218,10){\makebox(25,0.5){(1.1)}}
\put(25,20){\makebox(3,1){\scriptsize$i_\ast$}}
\put(25,10){\makebox(3,1){\scriptsize$i^!$}}
\put(74,20){\makebox(3,1){\scriptsize$j^\ast$}}
\put(74,10){\makebox(3,1){\scriptsize$j_\ast$}}
\end{picture}
\end{center}

\noindent  of functors between abelian categories such that

${\rm (i)}$ \ $i_*$ and $j^*$ are exact functors;

${\rm (ii)}$ \ $i_*$ and $j_*$  are fully faithful;

${\rm (iii)}$ \   $(i_*, i^!)$ and $(j^*, j_*)$ are adjoint pairs;
and

${\rm (iv)}$ \  ${\rm Im} i_* = {\rm Ker} j^*$ (and thus $i^!j_* =
0$).

\vskip5pt

In a right recollement $(1.1)$ the functor $i^!$ and $j_*$ are left
exact but not exact, in general. A right recollement is also called
{\it a localization sequence} e.g. in  [S], [G], [IKM], and [Kr],
and {\it a step} in [BGS].

\vskip5pt

{\it A left recollement} $(\mathcal B, \ \mathcal A, \ \mathcal C, \
i^*, \ i_*,  \ j_!, \ j^*)$ of abelian categories is a diagram

\begin{center}
\begin{picture}(100,28)
\put(13,10){\makebox(-22,2) {$\mathcal B$}}
\put(10,6){\vector(1,0){30}} \put(40,16){\vector(-1,0){30}}

\put(37,10){\makebox(25,0.8) {$\mathcal A$}}
\put(60,6){\vector(1,0){30}} \put(90,16){\vector(-1,0){30}}
\put(88,10){\makebox(25,0.5){$\mathcal C$}}
\put(218,10){\makebox(25,0.5){(1.2)}}\put(25,20){\makebox(3,1){\scriptsize$i^*$}}
\put(25,10){\makebox(3,1){\scriptsize$i_*$}}
\put(74,20){\makebox(3,1){\scriptsize$j_!$}}
\put(74,10){\makebox(3,1){\scriptsize$j^*$}}
\end{picture}
\end{center}
\noindent of functors between abelian categories such that

${\rm (i)}$ \  \ $i_*$ and $j^*$ are exact;

${\rm (ii)}$ \ $i_*$ and $j_!$ are fully faithful;

${\rm (iii)}$ \    $(i^*, i_*)$ and $(j_!, j^*)$ are adjoint pairs;
and

${\rm (iv)}$ \  ${\rm Im} i_* = {\rm Ker}j^*$ (and thus $i^*j_! =
0$).

\vskip5pt

Note that in a left recollement $(1.2)$ the functor $i^*$ and $j_!$
are right exact but not exact, in general. Thus, given a right
recollement $(\mathcal B, \ \mathcal A, \ \mathcal C, \ i_*, \ i^!,
\ j^*,\ j_*)$, the data  $(\mathcal C, \ \mathcal A, \ \mathcal B, \
j^*, \ j_*, \ \ i_*, \ i^!)$ is not a left recollement in general (similar
remark for a left recollement).

\vskip5pt

A recollement is first introduced for triangulated categories by A.
Beilinson, J. Bernstein, and P. Deligne [BBD]. A recollement of
abelian categories appeared in [Ku] and [CPS]. A {\it recollement}
$(\mathcal B, \ \mathcal A, \ \mathcal C, \ i^*, \ i_*, \ i^!,  \
j_!, \ j^*, \ j_*)$ of abelian categories is a diagram
\begin{center}
\begin{picture}(100,40)
\put(13,20){\makebox(-22,2) {$\mathcal B$}}
\put(37,20){\makebox(25,0.8) {$\mathcal A$}}
\put(88,20){\makebox(25,0.5){$\mathcal C$}}
\put(218,20){\makebox(25,0.5){$(1.3)$}}
\put(10,20){\vector(1,0){30}} \put(60,20){\vector(1,0){30}}
\put(25,23){\makebox(3,1){\scriptsize$i_\ast$}}
\put(74,23){\makebox(3,1){\scriptsize$j^\ast$}}
\put(40,12){\vector(-1,0){30}} \put(90,12){\vector(-1,0){30}}
\put(25,15){\makebox(3,1){\scriptsize$i^{!}$}}
\put(74,15){\makebox(3,1){\scriptsize$j_{*}$}}
\put(40,28){\vector(-1,0){30}} \put(90,28){\vector(-1,0){30}}
\put(25,31){\makebox(3,1){\scriptsize$i^{*}$}}
\put(74,31){\makebox(3,1){\scriptsize$j_{!}$}}
\end{picture}
\end{center}
\vskip-10pt \noindent of functors between abelian categories such
that

${\rm (i)}$ \ \ $(i^*, i_*), (i_*, i^!)$, $(j_!, j^*)$ and $(j^*,
j_*)$ are adjoint pairs;

${\rm (ii)}$ \ \ $i_*, j_!$ and $j_*$ are fully faithful; and

${\rm (iii)}$\ \ ${\rm Im}i_* = {\rm Ker}j^*$.

\vskip5pt

Thus in a recollement $(1.3)$ the functors $i_*$ and $j^*$ are
exact. So $(1.3)$ is a recollement if and only if the upper two rows
is a left recollement and the lower two rows is a right recollement.
By V. Franjou and  T. Pirashvili [FP], recollements of abelian
categories have some different properties from recollement of
triangulated categories. For example, ${\rm Ker} i^* \ne {\rm
Im}j_!$ and ${\rm Ker} i^! \ne {\rm Im} j_*$ in general, and
Parshall-Scott's theorem on comparison between two recollements of
triangulated categories ([PS, Thm. 2.5]) does not hold in general.
See also [Ps, PV, GYZ].

\vskip5pt

In a recollement of abelian categories, if $i^*$ and $i^!$ are
exact, then $j_!$ and $j_*$ are also exact (see Prop. \ref{rightrec}
and  \ref{leftrec}). The following result describes recollements of
abelian categories with exact functors.

\vskip5pt

\begin{thm}\label{h>4} \ Given a recollement $(1.3)$ of
abelian categories.  If $i^*$ and $i^!$ are exact, then $i^*\cong
i^!$ and $j_!\cong j_*,$ and $\mathcal A\cong \mathcal B\oplus
\mathcal C$.
\end{thm}

Quite surprising, we have

\begin{thm}\label{type (10)} \ Assume that $\mathcal A$ is a Grothendieck category.
Then any left recollement $(1.2)$ of abelian categories can be extended to a recollement of abelian categories$.$ \end{thm}

As a consequence, a colocalizing subcategory of a Grothendieck category is localizing. We stress that a right
recollement of abelian categories does not necessarily extend to a
recollement, and that a localizing subcategory of a Grothendieck category is not necessarily colocalizing. See
Subsection 5.2. On the other hand, W. Geigle and H. Lenzing [GL,
Prop. 5.3] have proved that any Serre subcategory $\mathcal S$ of
the finitely generated module category of an Artin algebra is always
localizing and colocalizing.

\vskip5pt

\begin{thm}\label{type} \ The type of a Serre subcategory of a Grothendieck category $\mathcal
A$ is in the list $$(0, 0), \ (0, -1), \ (1, -1), \ (0, -2), \ (1,
-2), \ (2, -1), \ (+\infty, -\infty);$$ and for each $(m, -n)$ in
this list, there exists a Serre subcategory such that its type is
$(m, -n);$ and if a Serre subcategory $\mathcal S$ is of type
$(+\infty, -\infty),$ then $\mathcal A\cong \mathcal S\oplus
(\mathcal A/\mathcal S)$ as categories.
\end{thm}

The main tools for proving Theorem \ref{type} are the work of
strongly hereditary torsion pairs by H. Tachikawa and K. Ohtake [TO;
O], the work of localizing subcategories by Geigle-Lenzing [GL], and
the argument on right (left) recollements of abelian categories,
especially Theorems \ref{h>4} and \ref{type (10)}. This result could
also be reformulated in terms of {\it the height of a ladder} of a
Grothendieck category (see [BGS], [AHKLY], [ZZZZ]). \ Theorems
\ref{h>4}, \ref{type (10)} and \ref{type}  do not hold in
triangulated categories.

\section{\bf Preliminaries}

Throughout $\mathcal A$ is an abelian category. A subcategory means
a full subcategory closed under isomorphisms. We will use the properties of a Grothendieck category $\mathcal A$:
it is {\it well-powered} ([M]) in the sense that for each object
$A \in \mathcal A$, the class of the subobjects of $A$ forms a set;
$\mathcal A$ has coproducts and products, enough injective objects; and the canonical morphism from
a coproduct to the corresponding product is a monomorphism (see [F], [Mit]).

\subsection{Serre subcategories} For Serre subcategories and
quotient categories we refer to [G], [Pop], and [GL]. A subcategory $\mathcal S$ of $\mathcal A$ is {\it
a Serre subcategory} if $\mathcal S$ is closed under subobjects,
quotient objects, and extensions. If $\mathcal S$ is a Serre
subcategory of $\mathcal A$ with the inclusion functor $i: \mathcal
S \rightarrow \mathcal A$, then we have the quotient category
$\mathcal A/\mathcal S$ which is abelian,  and the quotient functor
$Q: \mathcal A \rightarrow \mathcal A/\mathcal S$ is exact with
$Qi = 0$, and $Q$ has the universal property in the sense that if $F: \mathcal A\rightarrow \mathcal C$ is an exact functor between abelian categories with
$F i = 0$, then there exists a unique exact functor $G:\mathcal A/\mathcal S\rightarrow \mathcal C$ such that $F = GQ$.

A Serre subcategory $\mathcal S$ is {\it localizing}, if the
quotient functor $Q: \mathcal A \rightarrow \mathcal A/\mathcal S$
has a right adjoint $j_*$. In this case, $j_*$ is fully faithful
([GL, Prop. 2.2]). Dually, a Serre subcategory $\mathcal S$ is {\it
colocalizing}, if $Q$ has a left adjoint $j_!$. In this case, $j_!$
is fully faithful (the dual of [GL, Prop. 2.2]).

\vskip5pt

\subsection{Exact sequences of abelian categories} A sequence $0\rightarrow \mathcal{B} \stackrel{i_{*}}
\longrightarrow \mathcal{A} \stackrel{j^*}\longrightarrow
\mathcal{C} \rightarrow 0$ of exact functors between abelian
categories is an {\it exact sequence}, provided that there exists a Serre subcategory $\mathcal S$ of an abelian category
$\mathcal A'$ such that there is a commutative diagram $$\xymatrix{0\ar[r] &
\mathcal{B}\ar[r]^-{i_{*}}\ar[d]^-{\cong} &
\mathcal{A}\ar[r]^-{j^*}\ar[d]^-{\cong} & \mathcal{C}\ar[d]^-{\cong}
\ar[r]& 0 \\ 0\ar[r] & \mathcal{S}\ar[r]^-{i}&
\mathcal{A'}\ar[r]^-{Q}&\mathcal{A'}/\mathcal{S} \ar[r]& 0.}$$

A sequence $0\rightarrow \mathcal{B} \stackrel{i_{*}}
\longrightarrow \mathcal{A} \stackrel{j^*}\longrightarrow
\mathcal{C} \rightarrow 0$ of exact functors between abelian
categories is an exact sequence if and only if $i_*$ is fully
faithful, $i_*\mathcal{B}$ is a Serre subcategory of $\mathcal{A}$,
$j^* i_* = 0$, and $j^*$  has also the universal property. In this
case, we have ${\rm Im}i_{*} = {\rm Ker}j^*$.

\subsection{Torsion pairs}  For torsion pairs in an abelian category we refer to [D],
[J], and [TO]. A pair $(\mathcal T, \mathcal F)$ of subcategories of
$\mathcal A$ is {\it a torsion pair} ([D]),  if ${\rm Hom}(T, F) =
0$ for $T\in \mathcal T$ and $F\in \mathcal F$, and for each object
$A\in \mathcal A$, there is an exact sequence $0 \rightarrow T
\rightarrow A \rightarrow  F \rightarrow 0$ with $T \in \mathcal T$
and $F \in \mathcal F$. In this case, the exact sequence is called
{\it the $t$-decomposition} of $A$ with respect to $(\mathcal T,
\mathcal F)$.  A subcategory $\mathcal T$ (resp. $\mathcal F$) is
{\it a torsion class} (resp. {\it a torsionfree class}) if there
exists a subcategory $\mathcal F$ (resp. $\mathcal T$) such that
$(\mathcal T, \mathcal F)$ is a torsion pair. If   $(\mathcal T,
\mathcal F)$ is a torsion pair, then $\mathcal F = \mathcal
T^{\perp_{0}}$ and $\mathcal T = \ ^{\perp_{0}}\mathcal F,$ where
$\mathcal T^{\perp_{0}}: = \{A\in\mathcal A \ | \ {\rm Hom}(T, A) =
0, \ \forall \ T\in\mathcal T\}$, and $^{\perp_{0}}\mathcal F$ is
dually defined. By S. E. Dickson [D, Thm.2.3], if $\mathcal A$ is a
well-powered abelian category with coproducts and products, then a
subcategory $\mathcal T$ (resp. $\mathcal F$) is a torsion class
(resp. a torsionfree class) if and only if $\mathcal T$ (resp.
$\mathcal F$) is closed under quotient objects, extensions, and
coproducts (resp. under subobjects, extensions, and products).

\vskip5pt

A subcategory $\mathcal B$ is {\it  weakly localizing}, provided
that for each object $A$ of $\mathcal A$, there exists an exact
sequence
$$0\rightarrow B_1 \longrightarrow A \longrightarrow C\longrightarrow B_2\rightarrow 0$$
with $B_1\in \mathcal B,  \ B_2\in \mathcal B,$ and $C\in \mathcal
B^{\perp_{\le 1}}: = \{A\in\mathcal A \ | \ {\rm Hom}(B, A) = 0 =
{\rm Ext}^1(B, A), \ \forall \ B\in\mathcal B\}.$ By W. Geigle and
H. Lenzing [GL, Prop. 2.2], a Serre subcategory $\mathcal S$ is
localizing if and only if it is weakly localizing.   Dually, a
subcategory $\mathcal B$ is {\it  weakly colocalizing}, provided
that for each object $A$ of $\mathcal A$, there exists an exact
sequence
$$0\rightarrow B_1 \longrightarrow C\longrightarrow A\longrightarrow B_2\rightarrow 0$$
with $B_1\in \mathcal B,  \ B_2\in \mathcal B,$ and $C\in \
^{\perp_{\le 1}}\mathcal B:= \{A\in\mathcal A \ | \ {\rm Hom}(A, B)
= 0 = {\rm Ext}^1(A, B), \ \forall \ B\in\mathcal B\}.$ By the dual
of [GL, Prop. 2.2], a Serre subcategory $\mathcal S$ is colocalizing
if and only if it is weakly colocalizing.

\vskip5pt

Following H. Tachikawa and K. Ohtake [TO], a torsion pair $(\mathcal
T, \mathcal F)$ is {\it hereditary} (resp. {\it cohereditary}), if
$\mathcal T$ (resp. $\mathcal F$) is closed under subobjects (resp.
quotient objects); and it is {\it strongly hereditary} (resp. {\it
strongly cohereditary}), if $\mathcal T$ (resp. $\mathcal F$) is
weakly localizing (resp. weakly colocalizing). Every strongly
hereditary (resp. strongly cohereditary) torsion pair is hereditary
(resp. cohereditary) (see [TO, Prop. 1.7$^*$, 1.7]). K. Ohtake [O,
Thm. 2.6, 1.6] has proved that if $\mathcal A$ has enough injective
objects (resp. enough projective objects), then every hereditary
(resp. cohereditary) torsion pair is also strongly hereditary (resp.
strongly cohereditary) (see also [TO, Thm. 1.8$^*$, 1.8]).

\subsection{} Let $F: \mathcal C\rightarrow \mathcal A$ be a fully faithful functor between abelian categories. We
say that  $F$ is {\it Giraud}  if $F$ has a left adjoint which is an
exact functor. Dually, $F$ is  {\it coGiraud}  if  $F$ has a right
adjoint which is exact. By [TO, Coroll. 3.8, 2.8],  $F$ is Giraud
(resp. coGiraud) if and only if $F\mathcal C$ is a Giraud
subcategory (a coGiraud subcategory) of $\mathcal A$ in the sense of
[TO].

\begin{lem}\label{giraud} \ Let $\mathcal A$ be
an abelian category.

\vskip10pt

${\rm (1)}$ \ Given a Giraud functor  $j_*: \mathcal C\rightarrow
\mathcal A$ with an exact left adjoint $j^*: \mathcal A\rightarrow
\mathcal C$, there exists a functor $i^!: \mathcal A \rightarrow
{\rm Ker}j^*$,  such that $({\rm Ker}j^*, \ \mathcal A, \ \mathcal C,
\ i, \ i^!, \ j^*, \ j_*)$ is a right recollement, where $i:
{\rm Ker}j^*\rightarrow \mathcal A$ is the inclusion functor.

\vskip5pt

${\rm (1')}$ \ Given a coGiraud functor  $j_!: \mathcal C\rightarrow
\mathcal A$ with exact right adjoint $j^*: \mathcal A\rightarrow
\mathcal C$, there exists a functor $i^*: \mathcal A \rightarrow
{\rm Ker}j^*$,  such that $({\rm Ker}j^*, \ \mathcal A, \ \mathcal C,
\ i^*, \ i,  \ j_!, \ j^*)$ is a left recollement, where $i:
{\rm Ker}j^*\rightarrow \mathcal A$ is the inclusion functor.

\vskip10pt

$(2)$ \ Let $0\rightarrow \mathcal{B} \stackrel{i_{*}}
\longrightarrow \mathcal{A} \stackrel{j^*}\longrightarrow
\mathcal{C} \rightarrow 0$ be an exact sequence of abelian
categories.  If $j^*$ has a right adjoint $j_*$, then $j_*$ is
fully faithful, and there exists a functor $i^!: \mathcal A
\rightarrow {\rm Ker}j^*$ such that $({\rm Ker}j^*, \ \mathcal A, \
\mathcal C, \ i, \ i^!, \ j^*, \ j_*)$ is a right recollement,
where $i: {\rm Ker}j^*\rightarrow \mathcal A$ is the inclusion
functor.

\vskip5pt

${\rm (2')}$ \ Let $0\rightarrow \mathcal{B} \stackrel{i_{*}}
\longrightarrow \mathcal{A} \stackrel{j^*}\longrightarrow
\mathcal{C} \rightarrow 0$ be an exact sequence of abelian
categories. If $j^*$ has a left adjoint $j_!$, then $j_!$ is
fully faithful, and there exists a functor $i^*: \mathcal A
\rightarrow {\rm Ker}j^*$ such that $({\rm Ker}j^*, \ \mathcal A, \
\mathcal C, \ i^*, \ i,  \ j_!, \ j^*)$ is a left recollement,
where $i: {\rm Ker}j^*\rightarrow \mathcal A$ is the inclusion
functor.

\vskip10pt

${\rm (3)}$ \ If $(\mathcal B, \mathcal C)$ is a strongly hereditary torsion pair in $\mathcal A$. Then
$\mathcal B$ is a Serre subcategory of $\mathcal A$ and there is a right recollement of abelian category
\begin{center}
\begin{picture}(100,25)
\put(13,10){\makebox(-22,2) {$\mathcal B$}}
\put(10,16){\vector(1,0){30}} \put(40,6){\vector(-1,0){30}}

\put(37,10){\makebox(25,0.8) {$\mathcal A$}}
\put(60,16){\vector(1,0){30}} \put(90,6){\vector(-1,0){30}}
\put(93,10){\makebox(25,0.5){$\mathcal A/\mathcal B$}}
\put(25,20){\makebox(3,1){\scriptsize$i$}}
\put(25,10){\makebox(3,1){\scriptsize$i^!$}}
\put(74,20){\makebox(3,1){\scriptsize$Q$}}
\put(74,10){\makebox(3,1){\scriptsize$j_\ast$}}
\end{picture}
\end{center}
with ${\rm Im}j_* = B^{\perp_{\le 1}}$,  where $i$ is the inclusion functor and $Q$ is the quotient functor.

\vskip5pt

${\rm (3')}$ \ If $(\mathcal B, \mathcal C)$ is a strongly cohereditary torsion pair in $\mathcal A$. Then
$\mathcal C$ is a Serre subcategory of $\mathcal A$ and there is a left recollement of abelian category
\begin{center}
\begin{picture}(100,25)
\put(13,10){\makebox(-22,2) {$\mathcal C$}}
\put(10,6){\vector(1,0){30}}
\put(40,16){\vector(-1,0){30}}

\put(37,10){\makebox(25,0.8) {$\mathcal A$}}
\put(60,6){\vector(1,0){30}}
\put(90,16){\vector(-1,0){30}}
\put(93,10){\makebox(25,0.5){$\mathcal A/\mathcal C$}}
\put(25,20){\makebox(3,1){\scriptsize$i^*$}}
\put(25,10){\makebox(3,1){\scriptsize$i$}}
\put(74,20){\makebox(3,1){\scriptsize$j_!$}}
\put(74,10){\makebox(3,1){\scriptsize$Q$}}
\end{picture}
\end{center}
with ${\rm Im}j_!\cong \ ^{\perp_{\le 1}}C$ as categories, where $i$ is the inclusion functor and $Q$ is the quotient functor.
\end{lem} \noindent {\bf Proof.}  We only prove $(1)$, (2) and $(3)$. The assertions $(1')$, (2') and $(3')$ can be dually proved.

\vskip5pt

$(1)$ \ By assumption $j^*$ is exact, thus ${\rm Ker}j^*$ is an
abelian category and the inclusion functor $i: {\rm
Ker}j^*\rightarrow \mathcal A$ is exact. We claim that $i$ admits a
right adjoint $i^!: \mathcal A \rightarrow {\rm Ker}j^*$. In fact,
for any $A \in \mathcal A$, there is an exact sequence $0\rightarrow
{\rm Ker}\zeta_A \rightarrow A \stackrel{\zeta_A}\rightarrow
j_*j^*A$, where $\zeta: {\rm Id}_\mathcal A\rightarrow j_*j^*$ is
the unit of the adjoint pair $(j^*, j_*)$. Put $i^!A: = {\rm
Ker}\zeta_A$. Since $j^*$ is exact and $j^*\zeta_A$ is an
isomorphism, ${\rm Ker}\zeta_A \in {\rm Ker}j^*$. Thus, $i^!:
\mathcal A \rightarrow {\rm Ker}j^*$ defines a functor. For $B \in
{\rm Ker}j^*$ and $A\in \mathcal A$, since ${\rm Hom}_{\mathcal
A}(B, j_*j^*A) \cong {\rm Hom}_{\mathcal C}(j^*B, j^*A) =0$, by
applying the left exact functor ${\rm Hom}_{\mathcal A}(B, -)$ to
the exact sequence we get an isomorphism ${\rm Hom}(B, {\rm
Ker}\zeta_A) \cong {\rm Hom}(iB, A)$. This proves the claim, and
hence $({\rm Ker}j^*, \ \mathcal A, \ \mathcal C, \ i, \ i^!, \
j^*,\ j_*)$ a right recollement.

\vskip5pt

$(2)$ \ Since $0\rightarrow \mathcal{B} \stackrel{i_{*}}
\longrightarrow \mathcal{A} \stackrel{j^*}\longrightarrow
\mathcal{C} \rightarrow 0$ is an exact sequence of abelian
categories, without loss of generality, one may regard $j^*$ just as
the quotient functor $Q: \mathcal A\rightarrow \mathcal
A/i_*\mathcal B$. By [GL, Prop. 2.2], $j_*$ is fully faithful. Thus
$j_*: \mathcal C\rightarrow \mathcal A$ is a Giraud functor with the
exact left adjoint $j^*: \mathcal A\rightarrow \mathcal C$, and
hence the assertion follows from ${\rm (1)}$.

\vskip5pt

$(3)$ \ By [TO, Prop. 1.7$^*$], $(\mathcal B, \mathcal C)$ is a hereditary torsion pair. Thus
$\mathcal B$ is a Serre subcategory. Since by assumption $\mathcal B$ is weakly localizing, by [GL, Prop. 2.2],
the quotient functor $Q: \mathcal A \rightarrow \mathcal A/\mathcal B$ admits a right adjoint
$j_*: \mathcal A/\mathcal B \rightarrow \mathcal A$ which is fully faithful, and ${\rm Im}j_* = B^{\perp_{\le 1}}$.
So $j_*$ is a Giraud functor with an exact left adjoint $Q$. By $(1)$, there exists a functor $i^!: \mathcal A \rightarrow
{\rm Ker}Q = \mathcal B$,  such that $(\mathcal B, \ \mathcal A, \ \mathcal A/\mathcal B,
\ i, \ i^!, \ Q, \ j_*)$ is a right recollement. $\s$

\section{\bf Recollements of
abelian categories with exact functors}

\subsection{} The following proposition gives the properties of a right recollement of abelian
categories we need. Some of them are well-known for recollements of
abelian categories (see [FP], [Ps], [PV]).

\begin{prop}\label{rightrec} \ Let $(1.1)$ be a right recollement of abelian
categories. Then

\vskip5pt

$(1)$ \ ${\rm Im}i_*$ is a weakly localizing subcategory.
Explicitly, for each object $A \in \mathcal A$, there is an exact
sequence $0\rightarrow i_*i^!A \stackrel{\omega_A} \longrightarrow A
\stackrel{\zeta_A}\longrightarrow j_*j^*A \longrightarrow i_*B
\rightarrow 0$ for some object $B \in \mathcal B$, with $j_*j^*A\in
({\rm Im}i_*)^{\perp_{\le 1}}$, where $\omega$ the counit and
$\zeta$ is the unit.

\vskip5pt

$(2)$ \ $0\rightarrow \mathcal{B} \stackrel{i_{*}} \rightarrow
\mathcal{A} \stackrel{j^*}\rightarrow \mathcal{C} \rightarrow 0$ is
an exact sequence of abelian categories.

\vskip5pt

$(3)$ \ ${\rm Ker}i^! = ({\rm Im}i_*)^{\perp_0};$  $({\rm Im}i_*,
{\rm Ker}i^!)$ is a strongly hereditary torsion pair in $\mathcal
A,$ and $0\rightarrow i_*i^!A \stackrel{\omega_A} \longrightarrow A
\rightarrow {\rm Coker}  \omega_A \rightarrow 0$ is the
$t$-decomposition of $A$.

\vskip5pt

${\rm (4)}$ \ The following are equivalent$:$

 ${\rm (i)}$ \ \ $i^!$ is exact$;$

 ${\rm (ii)}$ \  $i^!$ and $j_*$ are exact$;$

 ${\rm (iii)}$ \  $0\rightarrow \mathcal{C} \stackrel{j_{*}} \rightarrow
\mathcal{A} \stackrel{i^!}\rightarrow \mathcal{B} \rightarrow 0$ is an
exact sequence of abelian categories$;$

 ${\rm (iv)}$ \  the sequence $0\rightarrow i_*i^!A
\stackrel{\omega_A} \longrightarrow A
\stackrel{\zeta_A}\longrightarrow j_*j^*A \rightarrow 0$ is exact
for each object $A \in \mathcal{A};$

${\rm (v)}$ \ \ ${\rm Im}j_* = {\rm Ker}i^!$,   $({\rm Im}i_*, {\rm
Im}j_*)$ is a cohereditary torsion pair in $\mathcal A,$ and
$0\rightarrow i_*i^!A \stackrel{\omega_A} \longrightarrow A
\stackrel{\zeta_A}\longrightarrow j_*j^*A \rightarrow 0$ is the
$t$-decomposition of $A;$

${\rm (vi)}$ \  $({\rm Im}i_*, {\rm Im}j_*)$ is a hereditary and
cohereditary torsion pair in $\mathcal A;$

${\rm (vii)}$ \ $({\rm Im}i_*, {\rm Im}j_*)$ is a strongly
hereditary and strongly cohereditary torsion pair in $\mathcal A;$

${\rm (viii)}$ \ $({\rm Im}i_*, {\rm Im}j_*)$ is a strongly
cohereditary torsion pair in $\mathcal A$.\end{prop}

\noindent {\bf Proof.}  $(1)$ \  Applying the exact functor $j^*$ to
the exact sequence $0 \rightarrow {\rm Ker}\zeta_{A}\rightarrow A
\stackrel{\zeta_{A}}\longrightarrow j_* j^* A
 \rightarrow {\rm Coker}\zeta_{A}\rightarrow 0$,
we get an exact sequence $0 \rightarrow j^*{\rm
Ker}\zeta_{A}\rightarrow j^* A \stackrel{j^*
\zeta_{A}}\longrightarrow j^* j_* j^* A
 \rightarrow j^*{\rm Coker}\zeta_{A}\rightarrow 0.$
Since  $j_*$ fully faithful, $j^* \zeta_{A}$ is an isomorphism. So
$j^*{\rm Ker}\zeta_{A} = 0 = j^*{\rm Coker}\zeta_{A}$, and then
${\rm Ker}\zeta_{A} \cong i_* B'$ and ${\rm Coker}\zeta_{A}\cong i_*
B$ for some $B' \in \mathcal B$ and $B\in \mathcal B$. Applying the
left exact functor $i_*i^!$ to the exact sequence $0 \rightarrow i_*
B'\rightarrow A \stackrel{\zeta_{A}}\rightarrow j_* j^* A$, we get a
commutative diagram \vspace {-5pt} $$\xymatrix{
  i_*i^! i_* B'\ar[d]^-{\omega_{i_*B'}} \ar[r]^-{\cong}&  i_*i^! A\ar[d]^-{\omega_A}\\
  i_* B'\ar[r]& A.}$$
Since $i_*$ is fully faithful, $\omega_{i_*B'}$ is an isomorphism,
and hence  ${\rm Ker}\zeta_{A} \cong i_* B'\cong i_*i^! A.$ So we
get the desired exact sequence.

\vskip5pt

To see ${\rm Im}i_*$ is weakly localizing, by the exact sequence
just established, it suffices to prove $j_*j^*A\in ({\rm
Im}i_*)^{\perp_{\le 1}}$ for each object $A\in\mathcal A$. We only
need to show ${\rm Ext}^1(i_*B, j_*j^*A) = 0$ for $B \in \mathcal
B$. Let $0 \rightarrow j_*j^*A \stackrel a\longrightarrow X
\rightarrow i_*B \rightarrow 0$ be an exact sequence. Applying the
exact functor $j^*$ we get an isomorphism $j^*a: j^*j_*j^*A \cong
j^*X $. By the commutative diagram
$$\xymatrix{
  \Hom(j^*X, j^*A) \ar[d]_-{\cong} \ar[r]^-{(j^*a, -)}_-{\cong} & \Hom(j^*j_*j^*A, j^*A ) \ar[d]^-{\cong} \\
  \Hom(X, j_*j^*A) \ar[r]^-{(a, -)} & \Hom(j_*j^*A, j_*j^*A)   }$$
we see that $\Hom(a, j_*j^*A): \Hom(X, j_*j^*A)\rightarrow
\Hom(j_*j^*A, j_*j^*A)$ is an isomorphism, which implies that the
exact sequence $0 \rightarrow j_*j^*A \stackrel a\longrightarrow X
\rightarrow i_*B \rightarrow 0$ splits.

\vskip5pt

$(2)$ \  Since ${\rm Im} i_* = {\rm Ker} j^*$ and $j^*$ is exact,
$i_*\mathcal B$ is a Serre subcategory of $\mathcal A$. Since
$j^*i_* = 0$, by the universal property of the quotient functor $Q:
\mathcal A \rightarrow \mathcal A/i_*\mathcal B$ we get a unique
exact functor $F: \mathcal A/  i_*\mathcal B \rightarrow \mathcal C$
such that $F Q \cong j^*$. We claim that $Qj_*$ is a quasi-inverse
of $F$. In fact,  $FQj_* \cong j^*j_* \cong {\rm Id}_{\mathcal C}$;
on the other hand, for each object $A\in \mathcal A$, by $(1)$ there
is a functorial isomorphism $QA\cong Qj_*j^*A$, and hence for each
object $QA\in \mathcal A/i_*\mathcal B$ there are functorial
isomorphisms
$$Qj_* F(QA) \cong Qj_* F(Qj_*j^*A)\cong Qj_*
(j^*j_*)j^*A\cong Qj_*j^*A\cong QA$$ in $\mathcal A/\mathcal B$,
i.e., $Qj_*F \cong {\rm Id}_{\mathcal A/i_*\mathcal B}$.

\vskip5pt

(3) \ For $A\in ({\rm Im}i_*)^{\perp_0}$, we have $\Hom_{\mathcal
B}(i^!A, i^!A) \cong \Hom_{\mathcal B}(i_*i^!A, A) =0$,  thus $i^!A
=0$. So $({\rm Im}i_*)^{\perp_0} \subseteq {\rm Ker}i^!.$ Conversely, if
$A\in {\rm Ker}i^!$, then for each $B\in\mathcal B$ we have
$\Hom_{\mathcal A}(i_*B, A) \cong \Hom_{\mathcal B}(B, i^!A) =0.$ So
${\rm Ker}i^! \subseteq ({\rm Im}i_*)^{\perp_0}.$ This shows ${\rm
Ker}i^! = ({\rm Im}i_*)^{\perp_0}$.

\vskip5pt

For $A\in\mathcal A,$ considering the exact sequence $0\rightarrow
i_*i^!A \stackrel{\omega_A} \longrightarrow A
\stackrel{\zeta_A}\longrightarrow j_*j^*A$ given in (1),  we get an
exact sequence $0\rightarrow i_*i^!A \stackrel{\omega_A}
\longrightarrow A \rightarrow {\rm Coker}  \omega_A  \rightarrow 0$.
To see $({\rm Im}i_*, {\rm Ker}i^!)$ is a torsion pair in $\mathcal
A,$ it suffices to show ${\rm Coker}  \omega_A  \in {\rm Ker}i^!.$
We see this by applying the left exact functor $i^!$ to the exact
sequence $0 \rightarrow {\rm Coker}  \omega_A \rightarrow j_*j^*A$,
and using $i^!j_* = 0$. By (1),  ${\rm Im}i_*$ is a weakly
localizing subcategory. Thus $({\rm Im}i_*, {\rm Ker}i^!)$ is a
strongly hereditary torsion pair.

\vskip5pt

(4) \ ${\rm (i)} \Rightarrow {\rm (ii)}:$  Applying the left exact
functor $j_*$ to a given exact sequence $0\rightarrow C_1
\stackrel{f} \longrightarrow C \stackrel{g}\longrightarrow C_2
\rightarrow 0$ in $\mathcal C$, we get an exact sequence
$0\rightarrow j_*C_1 \stackrel{j_* f} \longrightarrow j_* C
\stackrel{j_*g}\longrightarrow j_* C_2 \rightarrow {\rm Coker}(j_*
g) \rightarrow 0$; then by applying the exact functor $i^!$  we see
 $i^!{\rm Coker}(j_* g) = 0$ (since $i^!j_* C_2= 0$).  Applying the exact functor
$j^*$ we get an exact sequence
$$0\rightarrow j^*j_*C_1 \cong C_1 \rightarrow j^*j_* C \cong C
\rightarrow j^*j_* C_2 \cong C_2 \rightarrow j^* {\rm Coker}(j_* g)
\rightarrow0,$$ and thus $j^* {\rm Coker}(j_* g) = 0$. So ${\rm
Coker}(j_* g) = i_* B$ for some $B \in \mathcal B$, and hence $0 =
i^!{\rm Coker}(j_* g) = i^!i_*B \cong B$. Thus ${\rm Coker}(j_* g) =
i_* B = 0$, which proves the exactness of $j_*$.

\vskip5pt

${\rm (ii)}\Rightarrow {\rm (iii)}:$ We first claim ${\rm Im}j_* =
{\rm Ker}i^!$. It is clear that ${\rm Im}j_* \subseteq {\rm
Ker}i^!$. For each object $A \in {\rm Ker}i^!$, by $(1)$ we have an
exact sequence $0\rightarrow A \rightarrow j_*j^*A \rightarrow i_*B
\rightarrow 0;$ applying the exact functor $i^!$ we see that
$i^!i_*B=0,$ and hence $B \cong i^!i_*B = 0$. So $A\cong j_*j^*A\in
{\rm Im}j_*$. This proves the claim. Thus ${\rm Im}j_* = {\rm
Ker}i^!$ is a Serre subcategory. It remains to prove that $i^!$ has
the universal property. For this, assume that $F: \mathcal A
\rightarrow \mathcal B'$ is an exact functor such that $Fj_* = 0$.
Applying $F$ to the exact sequence in $(1)$ we get a functorial
isomorphism $Fi_*i^!(A) \cong F(A)$ for each object $A\in\mathcal
A$, i.e., $(Fi_*)i^!\cong F$. If $G: \mathcal B \rightarrow \mathcal
B'$ is an exact functor such that $Gi^!\cong F$, then $Gi^!\cong
(Fi_*)i^!$ and hence $G\cong Fi_*$ since $i^!$ is dense.

\vskip5pt

${\rm (iii)}\Rightarrow {\rm (iv)}:$  \ For each object $A \in
\mathcal A$,  applying the exact functor $i^!$ to the exact sequence
in $(1)$,  we get an exact sequence $0\rightarrow i^!i_*i^!A
\rightarrow i^!A\rightarrow i^!j_*j^*A \rightarrow i^!i_*B
\rightarrow 0.$ Since $i^!j_*j^*A = 0$,  $i^!i_*B = 0$. Thus $B
\cong i^!i_*B = 0$ and hence we get the desired exact sequence.

\vskip5pt

${\rm (iv)}\Rightarrow{\rm (v)}:$ From the given exact sequence one
easily see ${\rm Im}j_* = {\rm Ker}i^!$, and hence $({\rm Im}i_*,
{\rm Im}j_*)$ is a torsion pair by (3). It remains to prove that
${\rm Im}j_*= {\rm Ker}i^!$ is closed under quotient objects. For
this, let $0\rightarrow A_1 \rightarrow A\longrightarrow A_2
\rightarrow 0$ be an exact sequence with $A\in {\rm Ker}i^!$. Then
we get a commutative diagram with exact rows and columns:
$$\xymatrix@R=0.45cm{&&& 0\ar[d] & 0\ar[d] & 0\ar[d] & \\&&0\ar[r] & i_*i^!A_1\ar[r]^-{\omega_{A_1}}\ar[d] & A_1 \ar[r]^-{\zeta_{A_1}}\ar[d] & j_*j^*A_1\ar[r]\ar[d]
& 0 \\
&& 0\ar[r] & i_*i^!A\ar[r]^-{\omega_A}\ar[d] & A
\ar[r]^-{\zeta_A}\ar[d]^-g & j_*j^*A\ar[r]\ar[d]
& 0 &&&(3.1)\\
&& 0\ar[r] & i_*i^!A_2\ar[r]^-{\omega_{A_2}} & A_2
\ar[r]^-{\zeta_{A_2}}\ar[d] &
j_*j^*A_2\ar[r] & 0 \\
&&& & 0 &&}$$ Applying the Snake Lemma to the two columns on the
right, we get an exact sequence $0\rightarrow i_*i^!A_1
\longrightarrow i_*i^!A \longrightarrow i_*i^!A_2 \rightarrow 0,$
and hence $i_*i^!A_2 = 0.$ Since $i_*$ is  fully faithful, $i^!A_2 =
0$.

\vskip5pt

${\rm (v)}\Rightarrow{\rm (vi)}:$ Since ${\rm Im}i_* = {\rm Ker}
j^*$ and $j^*$ is exact,  ${\rm Im}i_*$ is closed under subobjects. So $({\rm Im}i_*, {\rm Im}j_*)$ is a hereditary and
cohereditary torsion pair.

\vskip5pt

${\rm (vi)}\Rightarrow {\rm (vii)}$ follows from [TO, Thm. 4.1] (we stress that this step does not need that $\mathcal A$ has enough injective objects).

\vskip5pt

${\rm (vii)}\Rightarrow {\rm (viii)}$ is clear.

\vskip5pt

${\rm (viii)}\Rightarrow {\rm (i)}:$ For each object $A \in \mathcal
A$, by $(1)$ we have an exact sequence $0\rightarrow i_*i^!A
\stackrel{\omega_A} \longrightarrow A
\stackrel{\zeta_A}\longrightarrow j_*j^*A \rightarrow i_*B
\rightarrow 0$ for some $B \in \mathcal B$. By assumption $({\rm
Im}i_*, {\rm Im}j_*)$ is a strongly cohereditary torsion pair, it
follows from [TO, Prop. 1.7] that $({\rm Im}i_*, {\rm Im}j_*)$ is a
cohereditary torsion pair. Thus ${\rm Im}j_*$ is closed under
quotient objects. So $i_*B\in {\rm Im}i_*\cap {\rm Im}j_* = \{0\},$
and hence we get the exact sequence $0\rightarrow i_*i^!A
\stackrel{\omega_A} \longrightarrow A
\stackrel{\zeta_A}\longrightarrow j_*j^*A \rightarrow 0$.

\vskip5pt

Let $0\rightarrow A_1 \rightarrow A \stackrel{g}\longrightarrow A_2
\rightarrow 0$ be an exact sequence in $\mathcal A$. Then we get a
commutative diagram (3.1) with exact rows and columns. Applying the
Snake Lemma to the two columns on the right, we get an exact
sequence
 $$0\rightarrow i_*i^!A_1
\longrightarrow i_*i^!A \stackrel{i_*i^!g}\longrightarrow i_*i^!A_2
\rightarrow 0.$$ Applying the left exact functor $i^!$ to
$0\rightarrow A_1 \rightarrow A \stackrel{g}\longrightarrow A_2
\rightarrow 0$ we get the exact sequence
 $$0\rightarrow i^!A_1 \longrightarrow i^!A \stackrel{i^!g}\longrightarrow i^!A_2 \longrightarrow {\rm Coker}(i^!g) \rightarrow 0,$$
and hence we have the exact sequence $0\rightarrow i_*i^!A_1
\longrightarrow i_*i^!A \stackrel{i_*i^!g}\longrightarrow i_*i^!A_2
\longrightarrow i_*{\rm Coker}(i^!g) \rightarrow 0.$ Thus $i_*{\rm
Coker}(i^!g) = 0$. Since $i_*$ is
fully faithful, ${\rm Coker}(i^!g) = 0$. This proves the exactness of $i^!$. $\s$

\vskip5pt

\subsection{} We state the dual result on left recollements without
proofs.

\begin{prop}\label{leftrec} \ Let $(1.2)$ be a left recollement of
abelian categories. Then

\vskip5pt

$(1)$ \ ${\rm Im}i_*$ is a weakly colocalizing subcategory.
Explicitly, for each object $A \in \mathcal{A}$, there is an exact
sequence $0\rightarrow i_*B \longrightarrow j_!j^*A
\stackrel{\epsilon_A} \longrightarrow A
\stackrel{\eta_A}\longrightarrow i_*i^*A \rightarrow 0$ for some $B
\in \mathcal B$, with $j_!j^*A\in \ ^{\perp_{\le 1}}({\rm Im}i_*)$,
where $\epsilon$ the counit and $\eta$ the unit.

 \vskip5pt

$(2)$ \ $0\rightarrow \mathcal{B} \stackrel{i_{*}} \rightarrow
\mathcal{A} \stackrel{j^*}\rightarrow \mathcal{C} \rightarrow 0$ is
an exact sequence of abelian categories.

\vskip5pt

${\rm (3)}$ \  ${\rm Ker}i^* = \ ^{\perp_0}({\rm Im}i_*);$ and
$({\rm Ker}i^*, {\rm Im}i_*)$ is a strongly cohereditary torsion
pair, and $0\rightarrow {\rm Ker}\eta_A \longrightarrow A
\stackrel{\eta_A}\longrightarrow i_*i^*A \rightarrow 0$ is the
$t$-decomposition of $A$.

\vskip5pt

${\rm (4)}$ \ The following are equivalent$:$

 ${\rm (i)}$ \ \ $i^*$ is exact$;$

 ${\rm (ii)}$ \  $i^*$ and $j_!$ are exact$;$

 ${\rm (iii)}$ \  $0\rightarrow \mathcal{C} \stackrel{j_!} \rightarrow
\mathcal{A} \stackrel{i^*}\rightarrow \mathcal{B} \rightarrow 0$ is
an exact sequence of abelian categories$;$

 ${\rm (iv)}$ \  the sequence $0\rightarrow j_!j^*A
\stackrel{\epsilon_A} \longrightarrow A
\stackrel{\eta_A}\longrightarrow i_*i^*A \rightarrow 0$ is exact for
each object $A \in \mathcal{A};$

${\rm (v)}$ \ \  ${\rm Im}j_! = {\rm Ker}i^*;$  $({\rm Im}j_!, {\rm
Im}i_*)$ is a hereditary torsion pair of $\mathcal A,$ and
$0\rightarrow j_!j^*A \stackrel{\epsilon_A} \longrightarrow A
\stackrel{\eta_A}\longrightarrow i_*i^*A \rightarrow 0$ is the
$t$-decomposition of $A;$

${\rm (vi)}$ \ \ $({\rm Im}j_!, {\rm Im}i_*)$ is a hereditary and
cohereditary torsion pair of $\mathcal A;$

${\rm (vii)}$ \ $({\rm Im}j_!, {\rm Im}i_*)$ is a strongly
hereditary and strongly cohereditary torsion pair of $\mathcal A;$

${\rm (viii)}$ \ $({\rm Im}j_!, {\rm Im}i_*)$ is a strongly
hereditary torsion pair of $\mathcal A$.\end{prop}

\begin{rem} By Lemma {\rm\ref{giraud}}$(1)$  $($resp. Lemma {\rm\ref{giraud}}$(1')$$)$, there is a bijective
correspondence between right $($resp. left$)$ recollements and
Giraud functors $($resp. coGiraud functors$)$.

\vskip5pt

By Lemma {\rm\ref{giraud}}$(3)$ and Proposition
{\rm\ref{rightrec}}$(3)$   $($resp. Lemma {\rm\ref{giraud}}$(3')$
and Proposition {\rm\ref{leftrec}}$(3)$$)$, there is a bijective
correspondence between right $($resp. left$)$ recollements and
strongly hereditary  $($resp. strongly cohereditary$)$ torsion
pairs.
\end{rem}

\subsection{} To prove Theorem \ref{h>4} we use the following fact,
in which the first assertion is just [TO, Lemmas 4.2, 4.2$^*$]. For
the use of the second assertion, we include a proof.

\begin{lem}\label{TFT}  \
Let  $(\mathcal{U}, \mathcal{V})$  and $(\mathcal{V}, \mathcal{W})$
be torsion pairs in $\mathcal{A}$. Assume that $\mathcal U$ is
closed under subobjects and $\mathcal W$ is closed under quotient
objects. Then $\mathcal{U} = \mathcal{W}$.

For each object $A \in \mathcal{A}$, let $0\rightarrow
A_{\mathcal{U}}\rightarrow A\stackrel{g}\rightarrow
A^{\mathcal{V}}\rightarrow 0$ be the $t$-decomposition of $A$ with
respect to $(\mathcal{U}, \mathcal{V})$, and
 $0\rightarrow A_{\mathcal{V}}\rightarrow A\stackrel{h}\rightarrow A^{\mathcal{W}}\rightarrow 0$ the $t$-decomposition of $A$ with respect to
 $(\mathcal{V}, \mathcal{W})$. Then $A_{\mathcal{V}}\cong A^{\mathcal{V}}$ and $A^{\mathcal W}\cong A_{\mathcal{U}}$,  \ $A\cong A_{\mathcal{U}}\oplus
A^{\mathcal{V}}$, and $A\mapsto(A_{\mathcal{U}}, A^{\mathcal{V}})$
gives an equivalence $\mathcal{A}\cong \mathcal{U}\oplus\mathcal{V}$
of categories.
\end{lem}
\noindent{\bf Proof.} Consider the push-out of $g$ and $h$, we get a
commutative diagram with exact rows and columns
$$\xymatrix @C=1.0pc @ R = 1.0pc{
& 0 \ar[d] & 0\ar[d] &0\ar[d] \\
0\ar[r] & E\ar[r]\ar[d] & A_{\mathcal{V}} \ar[r]\ar[d]&C \ar[d]\ar[r]& 0\\
    0\ar[r] & A_{\mathcal{U}}\ar[r]\ar[d] & A \ar[r]^-{g}\ar[d]^-{h}& A^{\mathcal{V}} \ar[r]\ar[d]& 0\\
    0\ar[r]& D \ar[r]\ar[d]& A^{\mathcal{W}}\ar[r]\ar[d]& B\ar[r]\ar[d] &0\\
    & 0 &0 &0
 }$$
Since $\mathcal{V}$ and $\mathcal{W}$ are closed under quotient
objects, $B \in \mathcal{V}\cap \mathcal{W} = \{0\}$. Since
$\mathcal{U}$ and $\mathcal{V}$ are closed under subobjects, $E\in
\mathcal{U}\cap \mathcal{V} = \{0\}$.

\vskip5pt

Thus $A_{\mathcal{V}} \cong C\cong A^{\mathcal{V}}$,  and hence $g$
is a splitting epimorphism
 $A \cong A_{\mathcal{U}}\oplus A^{\mathcal{V}}$. Also, $A_{\mathcal{U}} \cong D\cong
A^{\mathcal W}$, and hence $h$ is a splitting epimorphism. Taking
$A\in \mathcal U$ and $A\in \mathcal W$, respectively, we see that
$\mathcal{U}=\mathcal{W}$.  It is routine to
 verify that $\mathcal{A}\rightarrow \mathcal{U}\oplus\mathcal{V}$ given by $A\mapsto(A_{\mathcal{U}}, A^{\mathcal{V}})$ is an equivalence of categories.
$\s$

\vskip10pt

\noindent {\bf Proof of Theorem \ref{h>4}}. By Proposition
\ref{leftrec}$(4){\rm (v)}$, $({\rm Im}j_!, {\rm Im}i_*)$ is a
torsion pair with $t$-decomposition $0\rightarrow j_!j^*A
\stackrel{\epsilon_A} \longrightarrow A
\stackrel{\eta_A}\longrightarrow i_*i^*A \rightarrow 0$ and ${\rm
Im}j_! = {\rm Ker}i^*$ is closed under subobjects. By Proposition
\ref{rightrec}$(4){\rm (v)}$,  $({\rm Im}i_*, {\rm Im}j_*)$ is a
torsion pair with $t$-decomposition $0\rightarrow i_*i^!A
\stackrel{\omega_A} \longrightarrow A
\stackrel{\zeta_A}\longrightarrow j_*j^*A \rightarrow 0$ and ${\rm
Im}j_* = {\rm Ker} i^!$ is closed under quotient objects. It follows
from Lemma \ref{TFT} that ${\rm Im}j_! = {\rm Im}j_*$ and
$\mathcal{A} \cong {\rm Im}i_* \oplus {\rm Im}j_! \cong \mathcal{B}
\oplus \mathcal{C}$. For $A \in \mathcal{A}$, by Lemma \ref{TFT},
$i_*i^*A \cong i_*i^!A$ and  $j_!j^*A \cong j_*j^*A$. Since $i_*$ is
fully faithful, $i^*A \cong i^!A$ and hence $i^* \cong i^!$. Since
$j^*$ is dense, it follows that  $j_* \cong j_!$. $\s$

\section{\bf Proof of Theorem \ref{type (10)}}

\begin{lem}\label{tor<-leftrec} Let $\mathcal A$ be a
Grothendieck category, and $(1.2)$ a left recollement
of abelian categories. Then $(i_*\mathcal B, (i_* \mathcal
B)^{\perp_{0}})$ is a torsion pair.
\end{lem}
\noindent {\bf Proof.} By Proposition 3.2(3), $({\rm Ker}i^*,
i_*\mathcal B)$ is a torsion pair. So the torsionfree class
$i_*\mathcal B$ is closed under subobjects and products ([D, Thm.
2.3]). Since $\mathcal A$ is a Grothendieck category, any coproduct
is a subobject of the corresponding product, $i_*\mathcal B$ is also
closed under coproducts. On the other hand, since $i_*\mathcal B =
{\rm Ker}j^*$ and $j^*$ is exact, $i_*\mathcal B$ is closed under
quotient objects and extensions. Thus by [D, Thm. 2.3] $i_*\mathcal
B$ is a torsion class, and hence $(i_*\mathcal B, (i_* \mathcal
B)^{\perp_{0}})$ is a torsion pair. $\s$

\vskip10pt

\noindent{\bf Proof of Theorem \ref{type (10)}.} Given a left
recollement $(1.2)$, by Lemma \ref{tor<-leftrec}, $(i_*\mathcal B,
(i_* \mathcal B)^{\perp_{0}})$ is a torsion pair. Since $i_*\mathcal
B = {\rm Ker}j^*$ and $j^*$ is exact, $i_*\mathcal B$ is closed
under subobjects, i.e., $(i_*\mathcal B, (i_* \mathcal
B)^{\perp_{0}})$ is a hereditary torsion pair. Since $\mathcal A$ is
a Grothendieck category,  $\mathcal A$ has enough injective objects.
By [TO, Thm. 1.8$^*$], $(i_*\mathcal B, (i_* \mathcal
B)^{\perp_{0}})$ is a strongly hereditary torsion pair. Applying
Lemma \ref{giraud}(3) we get a right recollement
\begin{center}
\begin{picture}(100,25)
\put(10,10){\makebox(-22,2) {$i_*\mathcal B$}}
\put(10,16){\vector(1,0){30}} \put(40,6){\vector(-1,0){30}}

\put(37,10){\makebox(25,0.8) {$\mathcal A$}}
\put(60,16){\vector(1,0){30}} \put(90,6){\vector(-1,0){30}}
\put(97,10){\makebox(25,0.5){$\mathcal A/i_*\mathcal {B}$}}
\put(25,20){\makebox(3,1){\scriptsize${i}$}}
\put(25,10){\makebox(3,1){\scriptsize$\widetilde{i^!}$}}
\put(74,20){\makebox(3,1){\scriptsize$Q$}}
\put(74,10){\makebox(3,1){\scriptsize$j$}}
\end{picture}
\end{center}
where $i$ is the inclusion functor and $Q$ is the quotient functor. By Proposition 3.2(2), $0\rightarrow
\mathcal{B} \stackrel {i_*} \rightarrow \mathcal{A} \stackrel
{j^*}\rightarrow \mathcal{C} \rightarrow 0$ is an exact
sequence of abelian categories. By the universal property of the functors $j^*$ and $Q$, we get a commutative diagram
$$\xymatrix@R=0.55cm{
  \mathcal B \ar[d]_-{\cong}^-F \ar[r]^-{i_*} & \mathcal A \ar[d]_{=} \ar[r]^-{j^*} & \mathcal{C} \ar[d]_-{\cong}^-G \\
  i_*\mathcal B \ar[r]^-{i} & \mathcal A \ar[r]^-{Q} & \mathcal A/i_*\mathcal B}$$
and hence we get a recollement
  \begin{center}
\begin{picture}(100,40)
\put(13,20){\makebox(-22,2) {$\mathcal B$}}
\put(37,20){\makebox(25,0.8) {$\mathcal A$}}
\put(88,20){\makebox(25,0.5){$\mathcal C$}}
\put(10,20){\vector(1,0){30}} \put(60,20){\vector(1,0){30}}
\put(25,23){\makebox(3,1){\scriptsize$i_\ast$}}
\put(74,23){\makebox(3,1){\scriptsize$j^\ast$}}
\put(40,11){\vector(-1,0){30}} \put(90,11){\vector(-1,0){30}}
\put(25,15){\makebox(3,1){\scriptsize$i^{!}$}}
\put(74,15){\makebox(3,1){\scriptsize$j_{*}$}}
\put(40,28){\vector(-1,0){30}} \put(90,28){\vector(-1,0){30}}
\put(25,32){\makebox(3,1){\scriptsize$i^{*}$}}
\put(74,32){\makebox(3,1){\scriptsize$j_{!}$}}
\end{picture}
\end{center}
with $i^! = F^{-1}\widetilde{i^!}$ \  and $j_* = jG$. $\s$

\begin{cor}\label{colocaltolocal} A colocalizing subcategory of a
Grothendieck category is localizing.
\end{cor}
\noindent {\bf Proof.} Let $\mathcal S$ be a colocalizing
subcategory of a Grothendieck category $\mathcal A$. That is, the
quotient functor $Q: \mathcal A \rightarrow \mathcal A/\mathcal S$
has a left adjoint, denoted by $j_!: \mathcal A/\mathcal
S\rightarrow \mathcal A.$ By the dual of [GL, Prop. 2.2], $j_!$ is
fully faithful. So $j_!$ is a coGiraud functor with exact right
adjoint $Q$. By Lemma \ref{giraud}(1') there exists a functor $i^*:
\mathcal A \rightarrow {\rm Ker}Q = \mathcal S$,  such that
$(\mathcal S, \ \mathcal A, \ \mathcal A/\mathcal S, \ i^*, \ i,  \
j_!, \ Q)$ is a left recollement, where $i: \mathcal S\rightarrow
\mathcal A$ is the inclusion functor. Then by Theorem \ref{type
(10)} this left recollement can be extended to be a recollement, so
$Q$ has a right adjoint, i.e., $\mathcal S$ is localizing.
 $\s$

\section {\bf Proof of Theorem \ref{type}}

\subsection{\bf Serre subcategories of
type $(0, 0)$}

For a ring $R$, let ${\rm Mod}R$ be the category of right
$R$-modules. If  $R$ is a right noetherian, let ${\rm mod}R$ be the
category of finitely generated right $R$-modules.

\begin{lem}\label{00} Let R be a right noetherian ring. Then ${\rm mod}R$ is a Serre subcategory of type $(0, 0)$.
\end{lem}
\noindent {\bf Proof.}  It is clear that ${\rm mod}R$ is a Serre
subcategory of ${\rm Mod}R$. Let $i:{\rm mod}R\rightarrow {\rm
Mod}R$ and $Q:{\rm Mod}R \rightarrow {\rm Mod} R/{\rm mod}R$ be the
inclusion functor and the quotient functor, respectively. Assume
that the type of ${\rm mod}R$ is not $(0, 0)$. Then there exist
either adjoint pairs $(i_1, i)$ and $(j_1, Q)$, or adjoint pairs
$(i, i_{-1})$ and $(Q, j_{-1})$.

\vskip5pt

In the first case, by Lemma \ref{giraud}$(2')$ we get a
left recollement
\begin{center}
\begin{picture}(145,35)
\put(13,20){\makebox(-22,2) {${\rm mod}R$}}
\put(57,20){\makebox(25,0.8) {${\rm Mod}R$}}
\put(145,20){\makebox(25,0.5){${\rm Mod}R/{\rm mod}R.$}}
\put(20,16){\vector(1,0){30}} \put(90,16){\vector(1,0){30}}
\put(35,19){\makebox(3,1){\scriptsize$i$}}
\put(104,19){\makebox(3,1){\scriptsize$Q$}}
\put(50,26){\vector(-1,0){30}} \put(120,25){\vector(-1,0){30}}
\put(104,29){\makebox(3,1){\scriptsize$j_1$}}
\put(35,29){\makebox(3,1){\scriptsize$i_1$}}
\end{picture}
\end{center}
\vspace{-10pt}By Proposition \ref{leftrec}$(3)$ we have a torsion
pair $({\rm Ker}i_1, \ {\rm mod}R)$ in ${\rm Mod}R$. Thus the
torsionfree class ${\rm mod}R$ is closed under products, which is
absurd.

\vskip5pt

The dual argument shows that the second case is also impossible. We
give a direct proof. For each $X\in {\rm Mod}R/{\rm mod}R$ we have
${\rm Hom}_R(M, j_{-1}X) \cong {\rm Hom}_{{\rm Mod}R/{\rm mod}R}(QM,
X) = 0$ for all $M\in {\rm mod}R$. So $j_{-1}X$ has no non-zero
finitely generated submodule. Thus $j_{-1}X = 0$. Since $j_{-1}$ is
fully faithful, ${\rm Mod}R/{\rm mod}R = 0$, i.e., ${\rm mod}R =
{\rm Mod}R$, which is absurd.  $\s$

\subsection{\bf Serre subcategories of type $(0, -1)$}

\begin{lem}\label{0-1} Let ${\rm Ab}_t$ be the category of the torsion abelian groups. Then  ${\rm
Ab}_t$ is a Serre subcategory of type $(0, -1)$.
\end{lem}
\noindent {\bf Proof.}  \ Let ${\rm Ab}_f$ be the category of the
abelian groups in which every non-zero element is of infinite order.
Then $({\rm Ab}_t, {\rm Ab}_f)$ is a torsion pair in ModZ. Let
$i:{\rm Ab}_t \rightarrow {\rm Mod}Z$ and $Q: {\rm Mod}Z\rightarrow
{\rm Mod}Z/{\rm Ab}_t$ be the inclusion functor and the quotient
functor, respectively. Since ${\rm Ab}_t$ is closed under
submodules, it follows from [TO, Thm. 1.8*] (also [O, Thm. 2.6])
that $({\rm Ab}_t, {\rm Ab}_f)$ is a strongly hereditary torsion
pair. By Lemma \ref{giraud}$(3)$ we get a
right recollement
\begin{center}
\begin{picture}(130,25)
\put(8,10){\makebox(-22,2) {${\rm Ab}_t$}}
\put(10,15){\vector(1,0){30}} \put(40,6){\vector(-1,0){30}}

\put(45,9){\makebox(25,0.8) {${\rm Mod}Z$}}
\put(76,15){\vector(1,0){30}} \put(106,6){\vector(-1,0){30}}
\put(123,9){\makebox(25,0.5){${\rm Mod}Z/{\rm Ab}_t$}}
\put(25,18){\makebox(3,1){\scriptsize$i$}}
\put(25,9){\makebox(3,1){\scriptsize$i_{-1}$}}
\put(90,18){\makebox(3,1){\scriptsize$Q$}}
\put(90,9){\makebox(3,1){\scriptsize$j_{-1}$}}
\end{picture}
\end{center}
Thus $({\rm Ab}_t, {\rm Ker}i_{-1})$ is a torsion pair, by
Proposition \ref{rightrec}$(3)$. Comparing with the torsion pair
$({\rm Ab}_t, {\rm Ab}_f)$  we get ${\rm Ab}_f = {\rm Ker}i_{-1}$.

\vskip5pt

Assume that the type of ${\rm Ab}_t$ is not $(0, -1)$. Then there
exist either adjoint pairs $(i_1, i)$ and $(j_1, Q)$, or adjoint
pairs $(i_{-1}, i_{-2})$ and $(j_{-1}, j_{-2})$.

\vskip5pt

In the first case, by Lemma \ref{giraud}$(2')$ we get a left
recollement $({\rm Ab}_t, {\rm Mod}Z, {\rm Mod}Z/{\rm Ab}_t, \ i_1,
\  i, \ j_1, \ Q)$, and hence $({\rm Ker}i_1, {\rm Ab}_t)$ is a
torsion pair, by Proposition \ref{leftrec}$(3)$. So the torsionfree
class ${\rm Ab}_t$ is closed under products, which is absurd.

\vskip5pt

In the second case, the functor  $i_{-1}$ is exact, and hence ${\rm
Ab}_f = {\rm Ker}i_{-1}$ is closed under quotient groups, which is
absurd. $\s$

\vskip5pt

{\bf Remark.} \ The above argument also shows that there is a right
recollement of abelian categories which can not be extended to a
recollement (cf. Theorem \ref{type (10)}), and that a localizing
subcategory is not necessarily colocalizing (cf. Corollary
\ref{colocaltolocal}).

\subsection{\bf Serre subcategories
of type $(0, -2)$  and $(1, -1)$}

\begin{lem}\label{0-2} Let $\mathcal A$ be a Grothendieck category. Assume that both $(\mathcal T, \mathcal G)$ and $(\mathcal G, \mathcal F)$
are hereditary torsion pairs in $\mathcal A$, such that $\mathcal T$
is not a torsionfree class. Then  $\mathcal T$ is a Serre
subcategory of type $(0, -2)$, and $\mathcal G$ is a Serre
subcategory of type $(1, -1)$.
\end{lem}
\noindent {\bf Proof.}  It is clear that $\mathcal T$ and $\mathcal
G$ are Serre subcategories. Since  $\mathcal A$ is a Grothendieck
category, $\mathcal A$ has enough injective objects. Since
$(\mathcal T, \mathcal G)$ is a hereditary torsion pair, it follows
from [TO, Thm. 1.8$^*$]  that $(\mathcal T, \mathcal G)$ is strongly
hereditary.  By Lemma \ref{giraud}(3)
 there is a right recollement
\begin{center}
\begin{picture}(110,25)
\put(12,10){\makebox(-22,2) {$\mathcal T$}}
\put(10,15){\vector(1,0){30}} \put(40,6){\vector(-1,0){30}}
\put(37,9){\makebox(25,0.8) {$\mathcal A$}}
\put(60,15){\vector(1,0){30}} \put(90,6){\vector(-1,0){30}}
\put(92,9){\makebox(25,0.5){$\mathcal A/\mathcal T$}}
\put(25,19){\makebox(3,1){\scriptsize$i_\mathcal T$}}
\put(25,9){\makebox(3,1){\scriptsize$i_{-1}$}}
\put(72,19){\makebox(3,1){\scriptsize$Q_\mathcal T$}}
\put(72,9){\makebox(3,1){\scriptsize$j_{-1}$}}
\end{picture}
\end{center}
with ${\rm Im}j_{-1}= \mathcal T^{\perp_{\le 1}}$, where $i_\mathcal T$ and $Q_\mathcal T$ are respectively the inclusion functor and the quotient functor. We claim
$\mathcal G = \mathcal T^{\perp_{\le 1}}.$ In fact, $T^{\perp_{\le 1}}
\subseteq T^{\perp_{0}} = \mathcal G$. For each object $G\in
\mathcal G$, since $\mathcal T$ is a weakly localizing subcategory,
by definition there exists an exact sequence
$$0\rightarrow T_1 \stackrel a\longrightarrow G \longrightarrow C\stackrel b\longrightarrow T_2\rightarrow 0$$
such that $T_1\in \mathcal T,  \ T_2\in \mathcal T,$ and $C\in
\mathcal T^{\perp_{\le 1}}.$ But $(\mathcal T, \mathcal G)$ is a
torsion pair, $a = 0$ and $T_1 = 0$.  Since $\mathcal G$ is closed
under quotient objects, ${\rm Im} b\in \mathcal G$. Since by
assumption $\mathcal T$ is closed under subobjects, ${\rm Im} b\in
\mathcal T\cap \mathcal G = \{0\}$. So $G\cong C\in \mathcal
T^{\perp_{\le 1}}.$ This proves the claim.

\vskip5pt

Since $(\mathcal T, \mathcal G)$ is a hereditary  and cohereditary torsion pair, it follows from [TO, Thm. 4.1] that
$(\mathcal T, \mathcal G)$ is a
strongly cohereditary torsion pair.  By Lemma \ref{giraud}(3')
there is a left recollement
\begin{center}
\begin{picture}(110,25)
\put(12,10){\makebox(-22,2) {$\mathcal G$}}
\put(40,16){\vector(-1,0){30}} \put(10,6){\vector(1,0){30}}
\put(37,9){\makebox(25,0.8) {$\mathcal A$}}
\put(90,16){\vector(-1,0){30}} \put(60,6){\vector(1,0){30}}
\put(92,9){\makebox(25,0.5){$\mathcal A/\mathcal G$}}
\put(25,19){\makebox(3,1){\scriptsize$i^*$}}
\put(25,9){\makebox(3,1){\scriptsize$i_\mathcal G$}}
\put(72,19){\makebox(3,1){\scriptsize$j_!$}}
\put(72,9){\makebox(3,1){\scriptsize$Q_\mathcal G$}}
\end{picture}
\end{center}
with ${\rm Im}j_! = \ ^{\perp_{\le 1}}\mathcal G$,
where $i_\mathcal G$ and $Q_\mathcal G$ are respectively the inclusion functor and the quotient functor. Since we have shown
$\mathcal G = \mathcal T^{\perp_{\le 1}},$ it follows that
$\mathcal T \subseteq \ ^{\perp_{\le 1}}\mathcal G \subseteq \ ^{\perp_0}\mathcal G = \mathcal T.$ Thus
$\mathcal T = \ ^{\perp_{\le 1}}\mathcal G$.

\vskip5pt

Put $\widetilde{j_!}$ to be the equivalence $\mathcal A/\mathcal G
\rightarrow j_!(\mathcal A/\mathcal G) = \mathcal T$, and
$\widetilde{j_{-1}}$ to be the equivalence $\mathcal A/\mathcal T
\rightarrow j_{-1}(\mathcal A/\mathcal T) = \mathcal G.$ We claim
the diagram of functors
\begin{center}
\begin{picture}(120,30)
\put(12,10){\makebox(-22,2) {$\mathcal T$}}
\put(10,17){\vector(1,0){30}} \put(40,3){\vector(-1,0){30}}
\put(37,9){\makebox(25,0.8) {$\mathcal A$}}
\put(60,17){\vector(1,0){30}} \put(90,3){\vector(-1,0){30}}
\put(92,9){\makebox(25,0.5){$\mathcal A/\mathcal T$}}
\put(222,9){\makebox(25,0.5){$(5.1)$}}
\put(25,22){\makebox(3,1){\scriptsize$i_\mathcal T$}}
\put(25,9){\makebox(3,1){\scriptsize$\widetilde{j_!}Q_\mathcal G$}}
\put(72,22){\makebox(3,1){\scriptsize$Q_\mathcal T$}}
\put(72,9){\makebox(3,1){\scriptsize$i_\mathcal
G\widetilde{j_{-1}}$}}
\end{picture}
\end{center}
is a right recollement. In fact, since $j_{-1} = i_\mathcal
G\widetilde{j_{-1}}$, $(Q_\mathcal T, i_\mathcal
G\widetilde{j_{-1}})$ is an adjoint pair. By Proposition
\ref{leftrec}(1), for each object $A\in\mathcal A$ there is an exact
sequence
$$0\rightarrow i_\mathcal G G \longrightarrow j_!Q_\mathcal G A
\longrightarrow A
\longrightarrow i_\mathcal G i^*A \rightarrow 0$$
for some $G\in\mathcal G$. Since $j_!Q_\mathcal G A\in \mathcal T$ and $\mathcal T$ is closed under subobjects,
$i_\mathcal G G\in \mathcal T\cap \mathcal G = \{0\}$.
Thus for $T\in\mathcal T$ we have
${\rm Hom}(T, j_!Q_\mathcal G A)\cong {\rm Hom} (T, A).$ This shows that $(i_\mathcal T, \widetilde{j_!}Q_\mathcal G)$ is an adjoint pair.
This justifies the claim.

\vskip5pt

Again by  [TO, Thm. 1.8$^*$], $(\mathcal G, \mathcal F)$ is a
strongly hereditary torsion pair.  By Lemma \ref{giraud}(3)
there is a right recollement
\begin{center}
\begin{picture}(110,20)
\put(12,10){\makebox(-22,2) {$\mathcal G$}}
\put(10,15){\vector(1,0){30}} \put(40,5){\vector(-1,0){30}}
\put(37,9){\makebox(25,0.8) {$\mathcal A$}}
\put(60,15){\vector(1,0){30}} \put(90,5){\vector(-1,0){30}}
\put(92,9){\makebox(25,0.5){$\mathcal A/\mathcal G.$}}
\put(25,19){\makebox(3,1){\scriptsize$i_\mathcal G$}}
\put(25,9){\makebox(3,1){\scriptsize$i_{-2}$}}
\put(72,19){\makebox(3,1){\scriptsize$Q_\mathcal G$}}
\put(72,9){\makebox(3,1){\scriptsize$j_{-2}$}}
\end{picture}
\end{center}
Rewrite this we get a diagram of functors
\begin{center}
\begin{picture}(90,25)
\put(7,9){\makebox(-22,2) {$\mathcal A/\mathcal G$}}
\put(40,15){\vector(-1,0){30}} \put(10,5){\vector(1,0){30}}
\put(37,9){\makebox(25,0.8) {$\mathcal A$}}
\put(90,15){\vector(-1,0){30}} \put(60,5){\vector(1,0){30}}
\put(87,9){\makebox(25,0.5){$\mathcal G.$}}
\put(25,19){\makebox(3,1){\scriptsize$Q_\mathcal G$}}
\put(25,9){\makebox(3,1){\scriptsize$j_{-2}$}}
\put(72,19){\makebox(3,1){\scriptsize$i_\mathcal G$}}
\put(72,9){\makebox(3,1){\scriptsize$i_{-2}$}}
\end{picture}
\end{center}
(note that this is {\bf not} a left recollement, since $i_{-2}$ and $j_{-2}$ are not exact). Hence we have a diagram of functors
\begin{center}
\begin{picture}(110,25)
\put(12,9){\makebox(-22,2) {$\mathcal T$}}
\put(40,17){\vector(-1,0){30}} \put(10,3){\vector(1,0){30}}
\put(64,3){\vector(1,0){30}} \put(37,9){\makebox(25,0.8) {$\mathcal
A$}} \put(93,17){\vector(-1,0){30}} \put(64,3){\vector(1,0){30}}
\put(100,9){\makebox(25,0.5){$\mathcal A/\mathcal T.$}}
\put(222,9){\makebox(25,0.5){$(5.2)$}}
\put(25,24){\makebox(3,1){\scriptsize$\widetilde{j_!}Q_\mathcal G$}}
\put(25,9){\makebox(3,1){\scriptsize$j_{-2}\widetilde{j_!}^{-1}$}}
\put(76,9){\makebox(3,1){\scriptsize$\widetilde{j_{-1}}^{-1}i_{-2}$}}
\put(77,24){\makebox(3,1){\scriptsize$i_\mathcal
G\widetilde{j_{-1}}$}}
\put(76,9){\makebox(3,1){\scriptsize$\widetilde{j_{-1}}^{-1}i_{-2}$}}
\end{picture}
\end{center}

\vskip5pt

Putting $(5.1)$ and $(5.2)$ together we get a diagram of functors
\begin{center}
\begin{picture}(120,35)
\put(12,10){\makebox(-22,2) {$\mathcal T$}}
\put(10,22){\vector(1,0){30}} \put(40,8){\vector(-1,0){30}}
\put(37,9){\makebox(25,0.8) {$\mathcal A$}}
\put(60,22){\vector(1,0){30}} \put(90,8){\vector(-1,0){30}}
\put(228,9){\makebox(25,0.5){$(5.3)$}}
\put(95,9){\makebox(25,0.5){$\mathcal A/\mathcal T$}}
\put(25,26){\makebox(3,1){\scriptsize$i_\mathcal T$}}
\put(25,14){\makebox(3,1){\scriptsize$\widetilde{j_!}Q_\mathcal G$}}
\put(72,26){\makebox(3,1){\scriptsize$Q_\mathcal T$}}
\put(72,14){\makebox(3,1){\scriptsize$i_\mathcal
G\widetilde{j_{-1}}$}} \put(10,-5){\vector(1,0){30}}
\put(64,-5){\vector(1,0){30}}
\put(25,1){\makebox(3,1){\scriptsize$j_{-2}\widetilde{j_!}^{-1}$}}
\put(76,1){\makebox(3,1){\scriptsize$\widetilde{j_{-1}}^{-1}i_{-2}$}}
\end{picture}
\end{center}

\vskip10pt

\noindent such that $(i_\mathcal T, \ \widetilde{j_!}Q_\mathcal G, \ j_{-2}\widetilde{j_!}^{-1})$ and $(Q_\mathcal T, \ i_\mathcal G\widetilde{j_{-1}}, \ \widetilde{j_{-1}}^{-1}i_{-2})$ are adjoint sequences.

\vskip5pt

Assume that the type of $\mathcal T$ is not $(0, -2)$. Then there
exist either adjoint pairs $(i_1, i_\mathcal T)$ and $(j_1, Q_\mathcal T)$, or adjoint
pairs $(j_{-2}\widetilde{j_!}^{-1}, i_{-3})$ and $(\widetilde{j_{-1}}^{-1}i_{-2}, j_{-3})$.

\vskip5pt

In the first case, by Lemma \ref{giraud}$(2')$ we get a
left recollement $(\mathcal T, \mathcal A, \mathcal A/\mathcal T,
i_1, i_\mathcal T, j_1, Q_\mathcal T)$, and hence $({\rm Ker}i_1, \mathcal T)$ is a
torsion pair, by Proposition \ref{leftrec}$(3)$. This contradicts the assumption that $\mathcal T$ is not a torsionfree
class.

\vskip5pt

In the second case, all the functors in (5.3) are exact, and hence (5.3) is a recollement $(\mathcal A/\mathcal T, \ \mathcal A, \ \mathcal T)$.
By Theorem \ref{h>4} we have $i_\mathcal T\cong j_{-2}\widetilde{j_!}^{-1}$ and $Q_\mathcal T\cong\widetilde{j_{-1}}^{-1}i_{-2}$, and hence
both $i_\mathcal T$ and $Q_\mathcal T$ have left adjoints. This goes to the first case.

\vskip5pt

Thus the type of $\mathcal T$ is $(0, 2)$. This also proves the type of $\mathcal G$ is $(1, -1)$. $\s$

\vskip10pt

\begin{exm} Let $K$ be a field,  $R: = \prod\limits_{i=1}^\infty K_i$ and $I: =
\bigoplus\limits_{i=1}^\infty K_i$ with each $K_i = K$. Then $R$ is
a commutative ring and $I$ is an idempotent ideal of $R$. Put
$\mathcal G: = \{M\in {\rm Mod}R \ | \ MI = 0\}$. Then $\mathcal G$
is a {\rm TTF}-class in ${\rm Mod}R$,  i.e., $\mathcal G$ is a
torsion and torsion-free class, since $\mathcal G$ is subobjects,
quotient objects, extensions, coproducts and products. So we have a
{\rm TTF}-triple $(\mathcal T, \mathcal G, \mathcal F)$.

\vskip5pt

It is clear that $\mathcal T = \{M\in {\rm Mod}R \ | \ MI = M\}$. In
fact, for any $R$-module $M_1$ with $M_1I = M_1$ and $M_2 \in
\mathcal G $,  we have ${\rm Hom}(M_1, M_2) = 0;$ and for any $M \in
{\rm Mod}R$, we have an exact sequence $0\rightarrow MI \rightarrow
M \rightarrow M/MI \rightarrow 0$ with $(MI)I = MI$ and
$M/MI\in\mathcal G$.

\vskip5pt

We claim that $\mathcal T$ is closed under subobjects. By {\rm [D,
Thm. 2.9]} this is equivalent to say that $\mathcal G$ is closed
under taking injective envelopes. Thus, it suffices to prove that
for any $M \in \mathcal G$, the injective envelope $E(M)$ of $M$
satisfies $E(M)I=0$. Otherwise, there is an $m\in E(M)$ with $mI
\neq 0$. Set $L: = \{b \in I \ | \ mb \neq 0 \}.$ Then
$L\ne\emptyset$. Choosing  $b \in L$ such that the number of nonzero
components is smallest. We may assume that each nonzero component of
$b$ is $1_k$, the identity of $K$. In fact, if the nonzero
components of $b$ are exactly $b_{i_1}, \cdots, b_{i_t}$, where
$b_i$ is the $i$-th component of $b$, then we use $bb'$ to replace
$b$, where the nonzero components of $b'$ are exactly $b_{i_1}^{-1},
\cdots, b_{i_t}^{-1}$ $($note that $mbb'\ne 0:$ otherwise $mb =
mbb'b = 0$$)$. By the choice of $b$ we know that $mbI$ is a nonzero
submodule of $E(M)$. Since $E(M)$ is an essential extension of $M$,
$mbI\cap M\ne 0.$ So there is a nonzero element $r \in I$ with
$mbr\ne 0$ and $mbr\in M$. Note that the support of $br\in I$ is
contained in the support of $b$ $($by definition the support of $b$
is the set of $i$, such that the $i$-th component of $b$ is not
zero$)$. By the choice of $b$, the support of $br\in I$ is just  the
support of $b$. Let $b'_{i_1}, \cdots, b'_{i_t}$ be the nonzero
components of $br$, and  $d\in I$ the element with the nonzero
components exactly ${b'_{i_1}}^{-1}, \cdots, {b'_{i_t}}^{-1}$. Then
we get the desired contradiction $0\ne mb = mbrd \in MI = 0$.  This
proves the claim.

\vskip5pt

Since $K_i \in \mathcal T$ but $R\notin \mathcal T$, $\mathcal T =
\{M\in {\rm Mod}R \ | \ MI = M\}$ is not closed under products. Thus
$\mathcal T$ is not a torsion-free class. By Lemma {\rm\ref{0-2}},
the type of $\mathcal T$ is $(0, -2)$ and the type of $\mathcal G$
is $(1, -1)$.
\end{exm}

\subsection{\bf Serre subcategories of type $(1, -2)$ and $(2, -1)$}

Let $R$ and $S$ be rings,  $_SM_R$ a non-zero $S$-$R$-bimodule, and
$\Lambda=\left(\begin{smallmatrix}R & 0 \\ M &
S\end{smallmatrix}\right)$ the triangular matrix ring. A right
$\Lambda$-module is identified with a triple $(X, Y, f)$, where $X$
is a right $R$-module, $Y$ a right $S$-module, and $f: Y\otimes_S M
\rightarrow X$ a right $R$-map;  and a left $\Lambda$-module is
identified with a triple $\left(\begin{smallmatrix}U \\
V\end{smallmatrix}\right)_g$, where $U$ is a left $R$-module, $V$ a
left $S$-module, and $g: M\otimes_R U\rightarrow V$ a left $S$-map
([ARS, p.71]). Put $e_1= \left(\begin{smallmatrix}1 & 0 \\ 0 &
1\end{smallmatrix}\right)$ and $e_2= \left(\begin{smallmatrix}0 & 0
\\ 0 & 1\end{smallmatrix}\right)$.
 It is well-known that there is
a ladder of abelian categories (see [CPS, Sect. 2], [PV, 2.10]; also
[H, 2.1], [AHKLY, Exam. 3.4])
\begin{center}
\begin{picture}(140,40)
\put(13,20){\makebox(-22,2) {$ {\rm Mod}R$}}
\put(57,20){\makebox(25,0.8) {${\rm Mod}\Lambda$}}
\put(125,20){\makebox(25,0.5){${\rm Mod}S$}}
\put(20,20){\vector(1,0){30}} \put(90,20){\vector(1,0){30}}
\put(35,23){\makebox(3,1){\scriptsize$i_0$}}
\put(104,23){\makebox(3,1){\scriptsize$j_0$}}
\put(50,12){\vector(-1,0){30}} \put(120,12){\vector(-1,0){30}}
\put(104,15){\makebox(3,1){\scriptsize$j_{-1}$}}
\put(104,7){\makebox(3,1){\scriptsize$j_{-2}$}}
\put(50,29){\vector(-1,0){30}} \put(120,29){\vector(-1,0){30}}
\put(104,32){\makebox(3,1){\scriptsize$j_1$}}
\put(35,32){\makebox(3,1){\scriptsize$i_1$}}
\put(35,15){\makebox(3,1){\scriptsize$i_{-1}$}}
\put(20,4){\vector(1,0){30}} \put(90,4){\vector(1,0){30}}
\put(35,7){\makebox(3,1){\scriptsize$i_{-2}$}}
\end{picture}
\end{center}
i.e., the upper three rows form a recollement of abelian categories,
and $(i_{-1}, i_{-2})$ and $(j_{-1}, j_{-2})$ are adjoint pairs,
where
\begin{align*} & i_1=-\otimes_\Lambda \Lambda/\Lambda e_2\Lambda =-\otimes_\Lambda \left(\begin{smallmatrix}R \\
0
\end{smallmatrix}\right),
\\ & i_0={\rm Hom}_{\Lambda/\Lambda e_2\Lambda}(\Lambda/\Lambda e_2\Lambda, -) = {\rm Hom}_{R}(R,
-)\cong -\otimes_RR,
\\ & i_{-1} = {\rm Hom}_{\Lambda}(\Lambda/\Lambda e_2\Lambda, -) = {\rm Hom}_{\Lambda}(e_1\Lambda, -) \cong
-\otimes_\Lambda \Lambda e_1,
\\  & i_{-2} =  {\rm Hom}_{\Lambda / \Lambda e_2\Lambda}(\Lambda e_1, -) = {\rm Hom}_{R}(\left(\begin{smallmatrix}R \\
M
\end{smallmatrix}\right), -),
\\ & j_1 = -\otimes_{e_2\Lambda e_2} e_2\Lambda = -\otimes_S (M,
S), \\ &j_0 = {\rm Hom}_{\Lambda}(e_2\Lambda,-) = {\rm
Hom}_{\Lambda}((M,S), -)\cong -\otimes_\Lambda \Lambda e_2,
\\ & j_{-1} = {\rm Hom}_{e_2\Lambda e_2}(\Lambda e_2, -) = {\rm Hom}_{S}(S, -) \cong
-\otimes_SS, \\ & j_{-2} = {\rm Hom}_{\Lambda}(S, -),\end{align*}
where the right $\Lambda$-module $R$ is given by
$r\left(\begin{smallmatrix}r' & 0 \\ m & s\end{smallmatrix}\right):
= rr'$ and the right $\Lambda$-module $S$ is given by
$s\left(\begin{smallmatrix}r & 0 \\ m & s'\end{smallmatrix}\right):
= ss'$. Note that ${\rm Mod}R$ is a Serre subcategory of ${\rm
Mod}\Lambda$ and $0\rightarrow {\rm Mod}R \stackrel {i_0}\rightarrow
{\rm Mod}\Lambda \stackrel {j_0}\rightarrow {\rm Mod}S \rightarrow
0$ is an exact sequence of abelian categories.

\vskip5pt

We claim that the type of ${\rm Mod}R$ is $(1, -2)$.

\vskip5pt

In fact, since $M\ne 0$,  $\left(\begin{smallmatrix}R \\
0
\end{smallmatrix}\right)$ is not flat as a left $\Lambda$-module,  $i_1=-\otimes_\Lambda \left(\begin{smallmatrix}R \\
0
\end{smallmatrix}\right)$ is not exact, and hence $i_1$ has no left
adjoint. Also, since $M\ne 0$, $S$ is not projective as a right
$\Lambda$-module. So $j_{-2} = {\rm Hom}_{\Lambda}(S, -)$ is not
exact, and hence $j_{-2}$ has no right adjoint. This proves the
claim. (We include another proof. If both $i_1$ and $j_1$ have left
adjoints, then $i_1$ and
$j_1$ are exact, and hence  $i_1\cong i_{-1}$ by Theorem \ref{h>4}, i.e., $-\otimes_\Lambda \left(\begin{smallmatrix}R \\
0
\end{smallmatrix}\right)\cong -\otimes_\Lambda \left(\begin{smallmatrix}R \\
M
\end{smallmatrix}\right)$. But this is not true, since $M\ne 0$. Similarly, if both $i_{-2}$ and $j_{-2}$ have right adjoints, then $i_{-2}$ and $j_{-2}$ are exact,
and hence  $i_0\cong i_{-2}$ by Theorem \ref{h>4}, i.e., ${\rm Hom}_{R}(R, -)\cong {\rm Hom}_{R}(\left(\begin{smallmatrix}R \\
M
\end{smallmatrix}\right), -)$, which is absurd.)

\vskip5pt

The argument above also shows that $0\rightarrow {\rm Mod}S
\stackrel {j_{-1}}\rightarrow {\rm Mod}\Lambda \stackrel
{i_{-1}}\rightarrow {\rm Mod}R \rightarrow 0$ is an exact sequence
of abelian categories, and that the type of ${\rm Mod}S$ is $(2,
-1)$.

\begin{rem}\label{strange} Consider $\Lambda = T_2(R): = \left(\begin{smallmatrix}R & 0 \\ R &
R\end{smallmatrix}\right)$, then $i_{-2}\cong j_1$ and hence we have
adjoint sequences
$$(i_1, \ i_0, \ i_{-1}, \ i_{-2}\cong j_1, \ j_0, \ j_{-1}, \ j_{-2})$$
such that $i_1-\otimes_\Lambda \left(\begin{smallmatrix}R \\
0
\end{smallmatrix}\right)$ has no left adjoint, and $j_{-2}= {\rm Hom}_{\Lambda}(R, -)$ has no right
adjoint. Graphically we have
\begin{center}
\begin{picture}(140,65)
\put(13,20){\makebox(-22,2) {$ {\rm Mod}R$}}
\put(57,20){\makebox(25,0.8) {${\rm Mod}\Lambda$}}
\put(125,20){\makebox(25,0.5){${\rm Mod}R$}}
\put(20,20){\vector(1,0){30}} \put(90,20){\vector(1,0){30}}
\put(35,23){\makebox(3,1){\scriptsize$i_0$}}
\put(104,23){\makebox(3,1){\scriptsize$j_0$}}
\put(50,12){\vector(-1,0){30}} \put(120,12){\vector(-1,0){30}}
\put(104,15){\makebox(3,1){\scriptsize$j_{-1}$}}
\put(104,5){\makebox(3,1){\scriptsize$j_{-2}$}}
\put(50,29){\vector(-1,0){30}} \put(120,29){\vector(-1,0){30}}
\put(104,32){\makebox(3,1){\scriptsize$j_1$}}
\put(35,32){\makebox(3,1){\scriptsize$i_1$}}
\put(35,15){\makebox(3,1){\scriptsize$i_{-1}$}}
\put(20,0){\vector(1,0){30}}

\put(90,0){\vector(1,0){30}}

\put(35,5){\makebox(3,1){\scriptsize$i_{-2}\cong j_1$}}

\put(50,-10){\vector(-1,0){30}}

\put(35,-5){\makebox(3,1){\scriptsize$j_{0}$}}
\put(20,-20){\vector(1,0){30}}
\put(35,-15){\makebox(3,1){\scriptsize$j_{-1}$}}
\put(50,-28){\vector(-1,0){30}}
\put(35,-25){\makebox(3,1){\scriptsize$j_{-2}$}}
\put(90,37){\vector(1,0){30}}
\put(104,40){\makebox(3,1){\scriptsize$i_{-1}$}}
\put(104,48){\makebox(3,1){\scriptsize$i_{0}$}}
\put(120,45){\vector(-1,0){30}} \put(90,53){\vector(1,0){30}}
\put(104,57){\makebox(3,1){\scriptsize$i_{1}$}}
\end{picture}
\end{center}
\end{rem}

\vskip40pt

\subsection{\bf Serre subcategories of type $(+\infty, -\infty)$}

Let $\mathcal S$ and $\mathcal T$ be abelian categories. Then as a
subcategory of $\mathcal S\oplus \mathcal T$, $\mathcal S$ is a
Serre subcategory of type $(+\infty, -\infty)$. In fact, it is clear
that $(p_1, i_1, p_1)$ and $(i_2, p_2, i_2)$ are adjoint sequences,
where $i_1$ and $i_2$ are embeddings, and $p_1$ and $p_2$ are
projections. On the other hand, if $\mathcal S$ is a  Serre
subcategory of $\mathcal A$ and the type of $\mathcal S$ is
$(+\infty, -\infty)$, then by Theorem \ref{h>4} it is easy to see
$\mathcal A\cong \mathcal S\oplus\mathcal (\mathcal A/\mathcal S)$.

\subsection{\bf Proof of Theorem \ref{type}.} Let
$\mathcal B$ be a Serre subcategory of type $(m, -n)$, where $m$ and
$n$ are in the set $\{+\infty, 0, 1, 2, \cdots\}$. Denote by $i:
\mathcal B \rightarrow \mathcal A$ the inclusion functor  and $Q:
\mathcal A\rightarrow \mathcal A/ \mathcal B$ the quotient functor.
Put $h: = m+n+1$.

\vskip5pt

{\bf Claim 1.} If $h\ge 5$, then $(m, -n) = (+\infty, -\infty)$.

\vskip5pt

Assume that there is a diagram of functors
\begin{center}
\begin{picture}(100,45)
\put(13,20){\makebox(-22,2) {$\mathcal B$}}
\put(37,20){\makebox(25,0.8) {$\mathcal A$}}
\put(92,20){\makebox(25,0.5){$\mathcal A/\mathcal B$}}
\put(10,20){\vector(1,0){30}} \put(60,20){\vector(1,0){30}}
\put(25,23){\makebox(3,1){\scriptsize$i_2$}}
\put(74,23){\makebox(3,1){\scriptsize$j_2$}}
\put(40,12){\vector(-1,0){30}} \put(90,12){\vector(-1,0){30}}
\put(25,15){\makebox(3,1){\scriptsize$i_1$}}
\put(74,15){\makebox(3,1){\scriptsize$j_1$}}

\put(10,4){\vector(1,0){30}} \put(60,4){\vector(1,0){30}}
\put(25,7){\makebox(3,1){\scriptsize$i$}}
\put(74,7){\makebox(3,1){\scriptsize$Q$}}

\put(40,28){\vector(-1,0){30}}
\put(90,28){\vector(-1,0){30}}
\put(25,31){\makebox(3,1){\scriptsize$i_3$}}
\put(74,31){\makebox(3,1){\scriptsize$j_3$}}

\put(10,36){\vector(1,0){30}}
\put(60,36){\vector(1,0){30}}

\put(25,40){\makebox(3,1){\scriptsize$i_4$}}

\put(74,40){\makebox(3,1){\scriptsize$j_4$}}
\end{picture}
\end{center}
such that $(i_4, i_3, i_2, i_1, i)$ and $(j_4, j_3, j_2, j_1, Q)$
are adjoint sequences. Then $i_1, i_2, i_3, j_1,
j_2, j_3$ are exact. Since a left adjoint of $Q$ is fully faithful
(the dual of [GL, Prop. 2.2]), $j_1$ is fully faithful. Thus the two
rows at the bottom form a left recollement. By Proposition
\ref{leftrec}$(4){\rm (iii)}$, $0\rightarrow \mathcal A/\mathcal B
\stackrel{j_1}\rightarrow \mathcal A \stackrel{i_1}\rightarrow
\mathcal B \rightarrow 0$ is an exact sequence of abelian categories
and ${\rm Ker} i_1 = {\rm Im}j_1$. It follows from Lemma
\ref{giraud}$(2')$ that $i_2$ is fully faithful. Thus
\begin{center}
\begin{picture}(90,35)
\put(8,20){\makebox(-22,2) {$\mathcal A/\mathcal B$}}
\put(37,20){\makebox(25,0.8) {$\mathcal A$}}
\put(88,20){\makebox(25,0.5){$\mathcal B$}}
\put(10,20){\vector(1,0){30}}
\put(60,20){\vector(1,0){30}}
\put(25,23){\makebox(3,1){\scriptsize$j_1$}}
\put(74,23){\makebox(3,1){\scriptsize$i_1$}}
\put(40,12){\vector(-1,0){30}}
\put(90,12){\vector(-1,0){30}}
\put(25,15){\makebox(3,1){\scriptsize$Q$}}
\put(74,15){\makebox(3,1){\scriptsize$i$}}
\put(40,28){\vector(-1,0){30}}
\put(90,28){\vector(-1,0){30}}
\put(25,31){\makebox(3,1){\scriptsize$j_2$}}
\put(74,31){\makebox(3,1){\scriptsize$i_2$}}
\end{picture}
\end{center}
\vskip-10pt \noindent is recollement such that $j_2$ and $Q$ are
exact. It follows from  Theorem \ref{h>4} that $j_2\cong Q$ and $i_2
\cong i$, and hence the type is $(+\infty, -\infty)$.

\vskip5pt

Assume that there is a diagram of functors
\begin{center}
\begin{picture}(100,45)
\put(13,20){\makebox(-22,2) {$\mathcal B$}}
\put(37,20){\makebox(25,0.8) {$\mathcal A$}}
\put(92,20){\makebox(25,0.5){$\mathcal A/\mathcal B$}}
\put(40,20){\vector(-1,0){30}} \put(90,20){\vector(-1,0){30}}
\put(25,23){\makebox(3,1){\scriptsize$i_1$}}
\put(74,23){\makebox(3,1){\scriptsize$j_1$}}
\put(10,12){\vector(1,0){30}} \put(60,12){\vector(1,0){30}}
\put(25,15){\makebox(3,1){\scriptsize$i$}}
\put(74,15){\makebox(3,1){\scriptsize$Q$}}

\put(40,4){\vector(-1,0){30}} \put(90,4){\vector(-1,0){30}}
\put(25,7){\makebox(3,1){\scriptsize$i_{-1}$}}
\put(74,7){\makebox(3,1){\scriptsize$j_{-1}$}}

\put(10,28){\vector(1,0){30}} \put(60,28){\vector(1,0){30}}
\put(25,31){\makebox(3,1){\scriptsize$i_2$}}
\put(74,31){\makebox(3,1){\scriptsize$j_2$}}

\put(40,36){\vector(-1,0){30}}
\put(90,36){\vector(-1,0){30}}

\put(25,40){\makebox(3,1){\scriptsize$i_3$}}

\put(74,40){\makebox(3,1){\scriptsize$j_3$}}
\end{picture}
\end{center}
such that $(i_3, i_2, i_1, i, i_{-1})$ and $(j_3, j_2, j_1, Q,
j_{-1})$ are adjoint sequences. Then the four functors $i_1, i_2,
j_1, j_2$ are exact. By the dual of [GL, Prop. 2.2], $j_1$ is fully faithful. By the same
argument as above we know that the type is $(+\infty, -\infty)$.

\vskip5pt

Assume that there is a diagram of functors
\begin{center}
\begin{picture}(100,45)
\put(13,20){\makebox(-22,2) {$\mathcal B$}}
\put(37,20){\makebox(25,0.8) {$\mathcal A$}}
\put(92,20){\makebox(25,0.5){$\mathcal A/\mathcal B$}}
\put(10,20){\vector(1,0){30}} \put(60,20){\vector(1,0){30}}
\put(25,23){\makebox(3,1){\scriptsize$i$}}
\put(74,23){\makebox(3,1){\scriptsize$Q$}}
\put(40,12){\vector(-1,0){30}} \put(90,12){\vector(-1,0){30}}
\put(25,15){\makebox(3,1){\scriptsize$i_{-1}$}}
\put(74,15){\makebox(3,1){\scriptsize$j_{-1}$}}

\put(10,4){\vector(1,0){30}} \put(60,4){\vector(1,0){30}}
\put(25,7){\makebox(3,1){\scriptsize$i_{-2}$}}
\put(74,7){\makebox(3,1){\scriptsize$j_{-2}$}}

\put(40,28){\vector(-1,0){30}} \put(90,28){\vector(-1,0){30}}
\put(25,31){\makebox(3,1){\scriptsize$i_1$}}
\put(74,31){\makebox(3,1){\scriptsize$j_1$}}

\put(10,36){\vector(1,0){30}} \put(60,36){\vector(1,0){30}}

\put(25,40){\makebox(3,1){\scriptsize$i_2$}}

\put(74,40){\makebox(3,1){\scriptsize$j_2$}}
\end{picture}
\end{center}
such that $(i_2, i_1, i, i_{-1}, i_{-2})$ and $(j_2, j_1, Q, j_{-1},
j_{-2})$ are adjoint sequences. Then $i_1, i_{-1}, j_1, j_{-1}$ are
exact. By [GL, Prop. 2.2] and its dual, $j_{-1}$ and $j_1$ are fully
faithful. Thus the three rows in the middle form a recollement such
that $i_1$ and $i_{-1}$ are exact. It follows from Theorem \ref{h>4}
that $i_1\cong i_{-1}$ and $j_1\cong j_{-1}$, and hence the type is
$(+\infty, -\infty)$.

\vskip5pt

Assume that there is a diagram of functors
\begin{center}
\begin{picture}(100,45)
\put(13,20){\makebox(-22,2) {$\mathcal B$}}
\put(37,20){\makebox(25,0.8) {$\mathcal A$}}
\put(92,20){\makebox(25,0.5){$\mathcal A/\mathcal B$}}
\put(40,20){\vector(-1,0){30}} \put(90,20){\vector(-1,0){30}}
\put(25,23){\makebox(3,1){\scriptsize$i_{-1}$}}
\put(74,23){\makebox(3,1){\scriptsize$j_{-1}$}}
\put(10,12){\vector(1,0){30}} \put(60,12){\vector(1,0){30}}
\put(25,15){\makebox(3,1){\scriptsize$i_{-2}$}}
\put(74,15){\makebox(3,1){\scriptsize$j_{-2}$}}

\put(40,4){\vector(-1,0){30}} \put(90,4){\vector(-1,0){30}}
\put(25,7){\makebox(3,1){\scriptsize$i_{-3}$}}
\put(74,7){\makebox(3,1){\scriptsize$j_{-3}$}}

\put(10,28){\vector(1,0){30}} \put(60,28){\vector(1,0){30}}
\put(25,31){\makebox(3,1){\scriptsize$i$}}
\put(74,31){\makebox(3,1){\scriptsize$Q$}}

\put(40,36){\vector(-1,0){30}} \put(90,36){\vector(-1,0){30}}

\put(25,40){\makebox(3,1){\scriptsize$i_1$}}

\put(74,40){\makebox(3,1){\scriptsize$j_1$}}
\end{picture}
\end{center}
such that $(i_1, i, i_{-1}, i_{-2}, i_{-3})$ and $(j_1, Q, j_{-1},
j_{-2}, j_{-3})$ are adjoint sequences. Then $i_{-1}, i_{-2},
j_{-1}, j_{-2}$ are exact, and $j_{-1}$ is fully faithful. So
\begin{center}
\begin{picture}(110,25)
\put(13,10){\makebox(-22,2) {$\mathcal B$}}
\put(10,16){\vector(1,0){30}} \put(40,6){\vector(-1,0){30}}

\put(37,10){\makebox(25,0.8) {$\mathcal A$}}
\put(60,16){\vector(1,0){30}} \put(90,6){\vector(-1,0){30}}
\put(92,10){\makebox(25,0.5){$\mathcal A/\mathcal B$}}
\put(25,20){\makebox(3,1){\scriptsize$i$}}
\put(25,10){\makebox(3,1){\scriptsize$i_{-1}$}}
\put(74,20){\makebox(3,1){\scriptsize$Q$}}
\put(74,10){\makebox(3,1){\scriptsize$j_{-1}$}}
\end{picture}
\end{center}
is a right recollement.  By Proposition \ref{rightrec}$(4){\rm
(iii)}$, $0\rightarrow \mathcal A/\mathcal B
\stackrel{j_{-1}}\rightarrow \mathcal A \stackrel{i_{-1}}\rightarrow
\mathcal B \rightarrow 0$ is an exact sequence of abelian categories
and ${\rm Ker} i_{-1} = {\rm Im}j_{-1}$. It follows from Lemma
\ref{giraud}$(2)$ that $i_{-2}$ is fully faithful. Thus
\begin{center}
\begin{picture}(100,35)
\put(8,20){\makebox(-22,2) {$\mathcal A/\mathcal B$}}
\put(37,20){\makebox(25,0.8) {$\mathcal A$}}
\put(88,20){\makebox(25,0.5){$\mathcal B$}}
\put(10,20){\vector(1,0){30}} \put(60,20){\vector(1,0){30}}
\put(25,23){\makebox(3,1){\scriptsize$j_{-1}$}}
\put(74,23){\makebox(3,1){\scriptsize$i_{-1}$}}
\put(40,12){\vector(-1,0){30}} \put(90,12){\vector(-1,0){30}}
\put(25,15){\makebox(3,1){\scriptsize$j_{-2}$}}
\put(74,15){\makebox(3,1){\scriptsize$i_{-2}$}}
\put(40,28){\vector(-1,0){30}} \put(90,28){\vector(-1,0){30}}
\put(25,31){\makebox(3,1){\scriptsize$Q$}}
\put(74,31){\makebox(3,1){\scriptsize$i$}}
\end{picture}
\end{center}
\vskip-10pt \noindent is recollement such that $Q$ and $j_{-2}$ are
exact. It follows from  Theorem \ref{h>4} that $Q\cong j_{-2}$ and
$i \cong i_{-2}$, and hence the type is $(+\infty, -\infty)$.

\vskip5pt

Assume that there is a diagram of functors
\begin{center}
\begin{picture}(110,45)
\put(13,20){\makebox(-22,2) {$\mathcal B$}}
\put(37,20){\makebox(25,0.8) {$\mathcal A$}}
\put(92,20){\makebox(25,0.5){$\mathcal A/\mathcal B$}}
\put(10,20){\vector(1,0){30}} \put(60,20){\vector(1,0){30}}
\put(25,23){\makebox(3,1){\scriptsize$i_{-2}$}}
\put(74,23){\makebox(3,1){\scriptsize$j_{-2}$}}
\put(40,12){\vector(-1,0){30}} \put(90,12){\vector(-1,0){30}}
\put(25,15){\makebox(3,1){\scriptsize$i_{-3}$}}
\put(74,15){\makebox(3,1){\scriptsize$j_{-3}$}}
\put(10,4){\vector(1,0){30}} \put(60,4){\vector(1,0){30}}
\put(25,7){\makebox(3,1){\scriptsize$i_{-4}$}}
\put(74,7){\makebox(3,1){\scriptsize$j_{-4}$}}

\put(40,28){\vector(-1,0){30}} \put(90,28){\vector(-1,0){30}}
\put(25,31){\makebox(3,1){\scriptsize$i_{-1}$}}
\put(74,31){\makebox(3,1){\scriptsize$j_{-1}$}}

\put(10,36){\vector(1,0){30}} \put(60,36){\vector(1,0){30}}

\put(25,40){\makebox(3,1){\scriptsize$i$}}

\put(74,40){\makebox(3,1){\scriptsize$Q$}}
\end{picture}
\end{center}
such that $(i, i_{-1}, i_{-2}, i_{-3}, i_{-4})$ and $(Q, j_{-1},
j_{-2}, j_{-3}, j_{-4})$ are adjoint sequences. Then the six
functors $ i_{-1}, i_{-2}, i_{-3}, j_{-1}, j_{-2}, j_{-3}$ are
exact. By the same argument as above we know that the type is
$(+\infty, -\infty)$.

\vskip5pt

Up to now we have proved {\bf Claim 1.} So, from now on we assume
that $h\le 4$, i.e., $m+n \le 3$. Then the type $(m, -n)$ of
$\mathcal B$ is in the list
$$(3, 0), \ (2,
-1), \ (1, -2), \ (0, -3), \ (2, 0), \ (1, -1),  \ (0, -2), \ (1,
0), \ (0, -1)  \ (0, 0).$$

\vskip5pt

Assume that there is a diagram of functors
\begin{center}
\begin{picture}(110,25)
\put(13,10){\makebox(-22,2) {$\mathcal B$}}
\put(10,6){\vector(1,0){30}} \put(40,16){\vector(-1,0){30}}

\put(37,10){\makebox(25,0.8) {$\mathcal A$}}
\put(60,6){\vector(1,0){30}} \put(90,16){\vector(-1,0){30}}
\put(92,10){\makebox(25,0.5){$\mathcal A/\mathcal B$}}
\put(25,20){\makebox(3,1){\scriptsize$i_1$}}
\put(25,10){\makebox(3,1){\scriptsize$i$}}
\put(74,20){\makebox(3,1){\scriptsize$j_1$}}
\put(74,10){\makebox(3,1){\scriptsize$Q$}}
\end{picture}
\end{center}
such that $(i_1, i)$ and $(j_1, Q)$ are adjoint pairs. Then $j_1$ is
fully faithful, by the dual of [GL, Prop. 2.2]. So it is a left
recollement, and hence by Theorem \ref{type (10)} it can be extended
to be recollement. This shows that the type of $\mathcal B$ is not
in the set $\{(3, 0), \ (2, 0), \ (1, 0)\}.$

\vskip5pt

{\bf Claim 2.} The type of $\mathcal B$ can not be $(0, -3)$.

\vskip5pt

Otherwise, there is a diagram of functors
\begin{center}
\begin{picture}(110,45)
\put(13,20){\makebox(-22,2) {$\mathcal B$}}
\put(37,20){\makebox(25,0.8) {$\mathcal A$}}
\put(92,20){\makebox(25,0.5){$\mathcal A/\mathcal B$}}
\put(10,17){\vector(1,0){30}} \put(60,17){\vector(1,0){30}}
\put(25,20){\makebox(3,1){\scriptsize$i_{-2}$}}
\put(74,20){\makebox(3,1){\scriptsize$j_{-2}$}}
\put(40,9){\vector(-1,0){30}} \put(90,9){\vector(-1,0){30}}
\put(25,12){\makebox(3,1){\scriptsize$i_{-3}$}}
\put(74,12){\makebox(3,1){\scriptsize$j_{-3}$}}

\put(40,25){\vector(-1,0){30}} \put(90,25){\vector(-1,0){30}}
\put(25,28){\makebox(3,1){\scriptsize$i_{-1}$}}
\put(74,28){\makebox(3,1){\scriptsize$j_{-1}$}}

\put(10,33){\vector(1,0){30}} \put(60,33){\vector(1,0){30}}

\put(25,37){\makebox(3,1){\scriptsize$i$}}

\put(74,37){\makebox(3,1){\scriptsize$Q$}}
\end{picture}
\end{center}
such that $(i, i_{-1}, i_{-2}, i_{-3})$ and $(Q, j_{-1}, j_{-2},
j_{-3})$ are adjoint sequences. Then $i_{-1}, i_{-2}, j_{-1},
j_{-2}$ are exact, and $j_{-1}$ is fully faithful. So, the upper two
rows form a right recollement. Thus, by Proposition
\ref{rightrec}(4)(iii), $0\rightarrow \mathcal A/\mathcal B
\stackrel{j_{-1}}\rightarrow \mathcal A \stackrel{i_{-1}}\rightarrow
\mathcal B \rightarrow 0$ is an exact sequence of abelian
categories, and hence $i_{-2}$ is fully faithful, by Lemma
\ref{giraud}(2). So the upper three rows form a recollement
$(\mathcal A/\mathcal B, \mathcal A, \mathcal B, Q, j_{-1}, j_{-2},
i, i_{-1}, i_{-2})$, and then $i \cong i_{-2}$ and $Q \cong j_{-2}$
by Theorem \ref{h>4}. Thus the type of $\mathcal B$ is $(+\infty,
-\infty)$, which is absurd.

\vskip5pt

It remains to prove that for each  $(m, -n)$ in the list
$$(+\infty, -\infty), \ (2, -1), \ (1, -2), \ (1, -1), \ (0, -2), \ (0, -1),  \ (0, 0)$$ there exists a Serre subcategory of
$\mathcal A$ such that its type is $(m, -n)$. This is true by
Subsections 5.1-5.5. This completes the proof of Theorem \ref{type}.

\vskip10pt

\noindent {\small Jian Feng, School of Mathematics, Shanghai Jiao Tong University,  Shanghai
200240, China\\
Pu Zhang,  \ \ School of Mathematics, Shanghai Jiao Tong University,
Shanghai 200240, China}
\end{document}